\documentclass[preprint,11pt]{elsarticle}
\usepackage{amsfonts}
\usepackage{amsmath}
\usepackage{mathrsfs}
\usepackage[applemac]{inputenc}
\usepackage{times, ulem, amsfonts}
\usepackage{mathpazo, amsmath, amssymb, amsthm,color}

\def\({\left(}
\def\){\right)}
\def\N{{\mathbb N}}

\theoremstyle{plain}

\numberwithin{equation}{section}
\newtheorem{theorem}{Theorem}[section]
\newtheorem{proposition}[theorem]{Proposition}
\newtheorem{Lemma}[theorem]{Lemma}

\newtheorem{definition}[theorem]{Definition}

\theoremstyle{remark}
\newtheorem{remark}[theorem]{Remark}

\newcommand{\R}{{\mathbb R}}

\newcommand\Z{{\mathbb{Z}}}

\begin{document}
\baselineskip=24pt
\bibliographystyle{plainmma}
\begin{center}
{\Large\bf Global well-posedness for the 3-D density-dependent
liquid crystal flows in the critical Besov spaces}\\[2ex]
\footnote{$^*$This work is supported by NSFC under grant
 number 11171116 and 11401180, and by the Fundamental Research Funds for the
Central Universities of China under the grant number 2012ZZ0072.}

\footnote{$^{**}$
Email Address: pingxiaozhai@163.com (XP ZHAI); \
 yshli@scut.edu.cn (YS LI); \
 yanwei.scut@yahoo.com.cn (W YAN).}
Xiaoping ZHAI$^1$, Yongsheng LI$^1$ \ and \ Wei YAN$^2$\\[2ex]
 $^1$Department of Mathematics,
 South China University of Technology,
\\[0.5ex]
 Guangzhou, Guangdong 510640, P. R. China\\[2ex]
$^2$College of Mathematics and Information Science,
    Henan Normal University
\\[0.5ex]
 Xinxiang, Henan 453007,  P. R. China
\end{center}

\bigskip
\centerline{\bf Abstract}
In this paper,  we prove the local well-posedness of 3-D density-dependent
liquid crystal flows with initial data in the critical Besov spaces, without assumptions of small density variation.
Furthermore, if the initial density is close enough to a positive constant and the critical  Besov norm of the liquid crystal orientation field
 and the horizontal components of the initial velocity field polynomially small compared with the critical  Besov norm to the veritcal component of the initial velocity field, then the  system has a unique global solution.



\bigskip
\leftskip 0 true cm
\rightskip 0 true cm

\section{Introduction and the main result}
In this paper, we consider the  following density-dependent
liquid crystal flows in $\R^3$:
\begin{eqnarray}\label{1.1}
\left\{\begin{aligned}
&\rho_t+\mathrm{\mathrm{div}}(\rho u)=0,\\
&(\rho u)_t+\mathrm{\mathrm{div}}(\rho u\otimes u)-\mu\Delta u+\nabla P=-\nabla\cdot (\nabla d \odot \nabla d),\label{1}\\
&\partial_{t} d-\nu\Delta d+u \cdot \nabla d=|\nabla d|^2d,\quad |d|=1,\label{2}\\
&\mathrm{\mathrm{div}}u=0,\\
&\rho|_{t=0}=\rho_0(x),\quad u|_{t=0}=u_0(x),\quad d|_{t=0}=d_0(x),\quad |d_0(x)| = 1\quad \mathrm{in } \quad{\mathbb R}^3,\label{3}
\end{aligned}\right.
\end{eqnarray}
where $\rho$ is the density and $u$ is the velocity field, $d$ stand for the unit vector field
which represents the macroscopic orientations and  $\mu$ an $\nu$  are two positive viscosity constants. The term $\nabla d \odot \nabla d$ in
the equation of conservation of momentum denotes the $3\times3$ matrix whose $(i ,j)$-th entry is given by
$\partial_id \cdot \partial_jd$ $(1 \le i, j \le 3)$.

When $d$ is a given constant unit vector, (\ref{1.1}) is reduced to the well-known inhomogeneous incompressible Navier-Stokes
system, which has been studied by many researchers.
When the viscous coefficient equals some positive constant, Lady$\check{z}$enskaja and
Solonnikov \cite{ladyzenskaja} first established the unique resolvability  in a bounded domain
with homogeneous Dirichlet boundary condition for $u$; similar result was obtained by Danchin \cite{danchin2004} in $\R^n$
with initial data in the almost critical Sobolev spaces; Simon \cite{simon} proved the global existence of
weak solutions. In general, the global existence of weak solutions with finite energy  with
variable viscosity was proved by Lions in \cite{lions} (see also the references therein, and the monograph
\cite{antontsev}). Yet the regularity and uniqueness of such weak solutions is a big open question in the field
of mathematical fluid mechanics, even in two space dimensions when the viscosity depends on the
density. Except under the assumptions:
$$\rho_0\in L^{\infty}(\mathbb{T}^2),\quad \inf_{c>0}\|\frac{\mu(\rho_0)}{c}-1\|_{L^{\infty}(\mathbb{T}^2)}\le \varepsilon,\quad \mathrm{and} \quad u_0\in H^1(\mathbb{T}^2),$$
Desjardins \cite{desjardins1997} proved that Lions weak solution $(\rho, u)$ satisfies $u \in L^{\infty}((0, T);H^1(\mathbb{T}^2))$ and $\rho\in L^\infty((0, T) \times \mathbb{T}^2)$ for any $T < \infty.$  Moreover, with additional assumption on the initial density, he
could also prove that $u \in L^2((0, \tau );H^2(\mathbb{T}^2))$ for some short time $\tau$. Gui and Zhang \cite{guiguilong2009} proved the global wellposedness of
inhomogeneous Navier-Stokes system with initial data satisfying $\|\rho_0-1\|_{H^{s+1}}$ being sufficiently small and $u_0 \in H^s(\R^2) \cap \dot{H}^{-\varepsilon}(\R^2)$
for some $s > 2$ and $0 < \varepsilon < 1$. However, the exact size of $\|\rho_0-1\|_{H^{s+1}}$ was not given in [16].
Very recently, Danchin and Mucha \cite{danchin2013} proved that: given initial density $\rho_0\in L^\infty(\Omega) $ with
a positive lower bound and initial velocity $u_0 \in H^2(\Omega)$ for some bounded smooth domain of
$\R^n$, the Navier-Stokes system  with constant viscosity has a unique local solution. Furthermore, with the
initial density being close enough to some positive constant, for any initial velocity in two space
dimensions, and sufficiently small velocity in three space dimensions, they also proved its global
wellposedness.
When the density $\rho$ is away from zero, In \cite{abidi2007+1}, Abidi proved in general space dimension $n$ that: if
$1 < p < 2n, 0 < \mu < \mu(\rho)$, given $a_0 \in B_{p,1}^{n/p}(\R^n)$
 and $u_0 \in B_{p,1}^{-1+n/p}(\R^n)$ inhomogeneous Navier-Stokes system has a global solution
provided that $\|a_0\|_{B_{p,1}^{n/p}}+\|u_0\|_{B_{p,1}^{-1+n/p}}\le c_0$
 for some sufficiently small $c_0$. Moreover, this solution is
unique if $1 < p \le n$. This result generalized the corresponding results in [8, 9] and was improved by
Abidi and Paicu in \cite{abidi2007} with  $a_0 \in B_{q,1}^{n/q}(\R^n)$
 and $u_0 \in B_{p,1}^{-1+n/p}(\R^n)$
 for $p, q$ satisfying some technical
assumptions. Abidi, Gui and Zhang removed the smallness condition for $a_0$ in \cite{abidi2012}. Notice that
the main feature of the density space is to be a multiplier on the velocity space and this allows to
define the nonlinear terms in the system . Recently, Danchin and Mucha \cite{danchin2012} proved a more
general wellposedness result of this system  with $\mu(\rho) = \mu > 0 $ by considering very rough densities in some
multiplier spaces on the Besov spaces $B_{p,1}^{-1+n/p}(\R^n)$
for $1 < p < 2n$, which in particular completes
the uniqueness result in \cite{abidi2007+1} for $p \in (n, 2n)$ in the case when $\mu(\rho) = \mu > 0.$
On the other hand, motivated by \cite{guiguilong2010,paicu2011,zhangting} concerning the global wellposedness of 3-D incompressible
anisotropic Navier-Stokes system with the third component of the initial velocity
field being large, Paicu and Zhang \cite{paicu2012} proved that: given $a_0 \in B_{q,1}^{n/q}(\R^3)$
and $u_0 = (u^h_0 , u^3_0) \in B_{p,1}^{-1+3/p}(\R^3)$ for
$1 < q \le p < 6$ and $1/q-1/p\le1/3,$
with $\mu(\rho) = \mu > 0$ has a unique global solution as long as
\begin{align*}
\left(\mu \|a_0\|_{B_{q,1}^{3/q}}+\|u_0^h\|_{\dot{B}_{p,1}^{-1+3/p}}
\right)
\exp\left\{C_0\|u_0^3\|_{\dot{B}_{p,1}^{-1+3/p}}/{\mu^2} \right\}
\le c_0\mu
\end{align*}
for some sufficiently small $c_0$. We emphasize that the proof in \cite{paicu2012} used in a fundamental way the
algebraical structure of the momentum equation, namely div u = 0.

When the  density $\rho$ is a constant, (\ref{1.1}) reduce to be the classical nematic liquid crystal system
which was first first introduced by Lin \cite{linfanghua1995} as a simplification of the Ericksen-Leslie model \cite{ericksen,leslie}. As it is difficult to deal with the high nonlinear term $ |\nabla d|^2d$ in the third equation of (\ref{1.1}), Lin and Liu \cite{linfanghua1995,linfanghua1996} proposed to study an approximate model of liquid crystals system by Ginzburg-Landau
 function, i.e., $ |\nabla d|^2d$  is replaced by $-f(d) =-\nabla F(d) =-\nabla \frac{1}{{(\varepsilon^2)}}(|d|^2-1)^2$ with $\varepsilon>0$, and the third equation of (\ref{1.1}) is replaced by
$\partial_td + (u \cdot\nabla)d = \mu(\Delta d-f(d)).$
In this situation, Lin and Liu \cite{linfanghua1995} established the existence of global weak solution in dimensions 2 and 3 and assumptions that $u_0\in L^2$ and $d_0\in H^1$. Existence and uniqueness of global classical solutions were also obtained by them in dimension2 provided $u_0\in H^1$ and $d_0\in H^2$, and provided the viscosity  $\nu$ is large in dimension3. The partial regularity of suitable weak solutions was studied in \cite{linfanghua1996}. However, as pointed out in \cite{linfanghua1995}, the vanishing limit of $\varepsilon\rightarrow 0$ is an open and challenging problem.
For the full system , in two independent papers \cite{hong} and \cite{linfanghua2010}, the authors respectively established the global existence of Leray-Hopf type weak solutions to liquid crystals system in dimension2, and proved that the solutions are smooth away from at most finitely many singular times which is analogous to that for the heat flows of harmonic maps (see \cite{chang}). The global existence of weak solutions in dimension 3 is still now an open problem. For strong solutions, the global existence of strong solutions with small initial data and the local existence of strong solutions with any initial data of system  have been studied by Li and Wang \cite{lixiaoli}, Lin and Ding \cite{linjunyu}, Wang \cite{wangchangyou}, Hineman and Wang \cite{hineman}, and Hao and Liu \cite{haoyihang}. Recently, the uniqueness of weak solutions to this system has been studied by Lin and Wang \cite{linfanghua2010+1}. Blow up criteria for local smooth solutions has been studied by \cite{huang, liuqiao}and the references therein.

Very recently, Wen and Ding \cite{wenhuanyao} obtained the local strong solutions to the problem (\ref{1.1}) with vacuum. Fan, Zhou and Nakamura \cite{fanjishan} extended it to be global in 2D in a bounded smooth domain.

In this paper, we suppose that the initial density verifies $\inf
\rho_0(x) > 0$ and thus, by
the maximum principle for the transport equation, we have $\inf
\rho(t, x) > 0$. We define
$a ={1}/{\rho}-1 $ and assume $\mu=\nu=1$ for convenience
which allows us to work with the following system:
\begin{eqnarray}\label{1.2}
\left\{\begin{aligned}
&\partial_ta+u\cdot\nabla a=0,\\
&\partial_tu+u \cdot \nabla u+(1+a)(\nabla P-\Delta u)=-(1+a)\nabla\cdot(\nabla d \odot \nabla d),\\
&\partial_{t} d-\Delta d+u \cdot \nabla d=|\nabla d|^2d,\quad |d|=1,\\
&\mathrm{\mathrm{div}}u=0,\\
&a|_{t=0}=a_0(x),\quad u|_{t=0}=u_0(x),\quad d|_{t=0}=d_0(x),\quad |d_0(x)| = 1\quad \mathrm{in } \quad{\mathbb R}^3.
\end{aligned}\right.
\end{eqnarray}

\newpage
The main result of the present paper is stated in the following theorem:
\begin{theorem}\label{zhuyaodingli}
Let $a_0\in B_{2,1}^{3/2}(\R^3)$, $(u_0,d_0)\in \dot{B}_{2,1}^{1/2}(\R^3)\times \dot{B}_{2,1}^{3/2}(\R^3)$ with $\mathrm{div} u_0=0$, and \begin{align}\label{xiajie}
1+a_0\ge \underline{b}
\end{align}
for some positive constant $\underline{b}$. Then there exists a positive time $T$ such that (\ref{1.2})
has a unique local-in-time solution $(a, u, d)$  with
\begin{align*}
&a\in C_b([0, T ]; B_{2,1}^{3/2}(\R^3), \hspace{0.5cm} u\in C_b([0, T ]; \dot{B}_{2,1}^{1/2}(\R^3)\cap L^1([0,T];\dot{B}_{2,1}^{5/2}(\R^3)),\nonumber\\
&d\in C_b([0, T ]; \dot{B}_{2,1}^{3/2}(\R^3)\cap L^1([0,T];\dot{B}_{2,1}^{7/2}(\R^3)).
\end{align*}
Moreover, There exist a positive constant  $C$ such that
for any $(a_0,u_0,d_0)$ verifying
\begin{align}\label{1.3}
C\|a_0\|_{B_{2,1}^{3/2}}(1+\|u_0^3\|_{\dot{B}_{2,1}^{1/2}}+(\|u_0^h\|_{\dot{B}_{2,1}^{1/2}}+\|d_0\|_{\dot{B}_{2,1}^{3/2}})^2)\le1,
\end{align}
and
\begin{align}\label{1.4}
C(\|u_0^h\|_{\dot{B}_{2,1}^{1/2}}+\|d_0\|_{\dot{B}_{2,1}^{3/2}})\left(1+(\|u_0^h\|_{\dot{B}_{2,1}^{1/2}}+\|d_0\|_{\dot{B}_{2,1}^{3/2}})^2+ \|u_0^3\|_{\dot{B}_{2,1}^{1/2}}^{1/2}(\|u_0^h\|_{\dot{B}_{2,1}^{1/2}}+\|d_0\|_{\dot{B}_{2,1}^{3/2}})^{1/2}\right)\le1.
\end{align}
Then the system $\mathrm{(\ref{1.1})}$ has a unique global solution $(a,u,d)$ with
$$a\in C([0,+\infty);B_{2,1}^{3/2} )\cap
\tilde{L}^\infty((0,+\infty);B_{2,1}^{3/2}), $$
$$u\in C([0,+\infty);\dot{B}_{2,1}^{1/2})\cap
\tilde{L}^\infty((0,+\infty);\dot{B}_{2,1}^{1/2})\cap
L^1((0,+\infty);\dot{B}_{2,1}^{5/2}),$$
$$d\in C([0,+\infty);\dot{B}_{2,1}^{3/2})\cap
\tilde{L}^\infty((0,+\infty);\dot{B}_{2,1}^{3/2})\cap
L^1((0,+\infty);\dot{B}_{2,1}^{7/2}).$$

\end{theorem}

The remainder of this paper is organized as follows. In section 2, we shall first review some useful statements on functional spaces and basic analysis tools, and introduce several technical Lemmas. In Section 3, we give the estimates of the transport equation. In Section 4, we  shall prove the local well-posedness of Theorem 1.1.  In last section, we complete the proof of the global wellposedness.

$\mathbf{Notations:}$ Let $A$, $B$ be two operators, we denote $[A; B] = AB - BA$, the commutator
between $A$ and $B$. For $a\lesssim b$, we mean that there is a uniform constant $C$, which may
be different on different lines, such that $a \le C b$.
For $X$ a Banach space and $I$ an interval of $\mathbb{R}$, we denote by $C(I; X)$ the set of
continuous functions on $I$ with values in $X$. For $q \in [1, +\infty]$, the notation $L^q (I; X)$ stands for the set of measurable
functions on $I$ with values in $X$, such that $t \rightarrow \|f(t)\|_{ X }$ belongs to $L^q (I)$. For a vector
$v = (v_1 , v_2 ) \in X$, we mean that all the components $v_i (i = 1, 2)$ of $v$ belong to the
space $X$.
We always denote $(d_j)_{j\in\mathbb{Z}}$ is a
generic element of $l^1(\mathbb{Z})$  so that $\sum_{j\in\mathbb{Z}}d_j=1$.

\section{Preliminaries }

Let  $(\chi,\phi)$  be two  smooth radial functions,
$0\leq  (\chi,\phi) \leq 1,$
such that $\chi$  is supported in the ball
$\mathcal{B}=\{\xi \in \R^{3}, |\xi|\leq {4}/{3}\}$ and $\varphi$  is
supported in the ring
 $\mathcal{C}=\{\xi \in \R^{3}, {3}/{4}\leq |\xi|\leq {8}/{3}\}$. Moreover, there holds
 \begin{align*}
 \sum_{j \in \mathbb{Z}} \varphi(2^{-j}\xi)=1, \hspace{0.5cm}\forall \xi \neq0.
 \end{align*}

Let $h=\mathcal {F}^{-1}\varphi$ and $\tilde{h}=\mathcal {F}^{-1}\chi$, then we define the dyadic blocks as follows:
\begin{align*}
&\dot{\Delta}_jf=\varphi(2^{-j}D)f=2^{3j}\int_{\R^3}h(2^{j}y)f(x-y)dy,\\
&\dot{S}_jf=\chi(2^{-j}D)f=2^{3j}\int_{\R^3}\tilde{h}(2^{j}y)f(x-y)dy.
\end{align*}
By telescoping the series, we thus have the following Littlewood-Paley decomposition
\begin{align*}
u=\
\sum_{j\in\mathbb{Z}}\dot{\Delta}_ju,\hspace{0.5cm} \forall u\in S'(\R^3)/\mathcal{P}[\R^3],
\end{align*}
where $\mathcal{P}[\R^3]$ is the set of polynomials (see\cite{bcd}).
Moreover, the Littlewood-Paley decomposition
satisfies the property of almost orthogonality:
$$
\dot{\Delta}_k\dot{\Delta}_ju\equiv 0 \hspace{0.5cm}\mathrm{if} \hspace{0.2cm}
|k-j|\ge2 \hspace{0.5cm}\mathrm{and}\hspace{0.5cm}
\dot{\Delta}_k(\dot{S}_{j-1}u\dot{\Delta}_ju)\equiv0\hspace{0.5cm} \mathrm{if}\hspace{0.2cm} |k-j|\ge5.
$$
 Now we recall the definition of nonhomogeneous Besov spaces.
\begin{definition}
Let $s\le 3/p$ (respectively $s \in {\mathbb R}$), $(r,\lambda,p) \in [1,+\infty]^3$ and $T \in (0,+\infty]$. We define
${\tilde{L}^\lambda_{T}(\dot{B}_{p,r}^s({\mathbb R}^3))}$ as the completion of $C([0,T ];\mathcal{S}({\mathbb R}^3))$ by the norm
$$
\|f\|_{\tilde{L}^\lambda_{T}(\dot{B}_{p,r}^s)}=\left\{\sum_{q\in {\mathbb Z}}2^{rqs}
\left(\int_0^T\|\dot{\Delta}_qf(t)\|_{L^p}^\lambda dt\right)^{ r/\lambda}\right\}^{ 1/r}<\infty,
$$
with the usual change if $r =\infty$. For short, we just denote this space by ${\tilde{L}^\lambda_{T}(\dot{B}_{p,r}^s)}$.
\end{definition}
\begin{remark}
It is easy to observe that for $p,r,\lambda,\lambda_1,\lambda_2\in[1,+\infty]$, $\theta\in[0,1]$,  $0<s_1<s_2,$ we have the following interpolation inequality in the Chemin-Lerner space (see\cite{bcd}):
\begin{align*}
\|u\|_{\tilde{L}^\lambda_{T}(\dot{B}_{p,r}^s)}\le\|u\|^\theta_{\tilde{L}^{\lambda_1}_{T}(\dot{B}_{p,r}^{s_1})}
\|u\|^{1-\theta}_{\tilde{L}^{\lambda_2}_{T}(\dot{B}_{p,r}^{s_2})}
\end{align*}
with ${1}/{\lambda}={\theta}/{\lambda_1}+{(1-\theta)}/{\lambda_2}$ and $s=\theta s_1+(1-\theta)s_2$.
\end{remark}
Let us emphasize that, according to the Minkowski inequality, we have
\begin{align*}
\|u\|_{\tilde{L}^\lambda_{T}(\dot{B}_{p,r}^s)}\le\|u\|_{L^\lambda_{T}(\dot{B}_{p,r}^s)}\hspace{0.2cm} \mathrm{if }\hspace{0.2cm}  \lambda\ge r,\hspace{0.2cm}
\|u\|_{\tilde{L}^\lambda_{T}(\dot{B}_{p,r}^s)}\ge\|u\|_{L^\lambda_{T}(\dot{B}_{p,r}^s)},\hspace{0.2cm} \mathrm{if }\hspace{0.2cm}  \lambda\le r.
\end{align*}
The following Bernstein's Lemma will be repeatedly used throughout this paper.
\begin{Lemma}\label{bernstein}
Let $\mathcal{B}$ be a ball and $\mathcal{C}$ a ring of $\R^3$. A constant $C$ exists so that for any positive real number $\lambda$, any
non-negative integer k, any smooth homogeneous function $\sigma$ of degree m, and any couple of real numbers $(a, b)$ with
$1\le a \le b$, there hold
\begin{align*}
&&\mathrm{Supp} \,\hat{u}\subset\lambda \mathcal{B}\Rightarrow\sup_{|\alpha|=k}
\|\partial^{\alpha}u\|_{L^b}\le C^{k+1}\lambda^{k+3(1/a-1/b)}\|u\|_{L^a},\\
&&\mathrm{Supp} \,\hat{u}\subset\lambda \mathcal{C}\Rightarrow C^{-k-1}\lambda^k\|u\|_{L^a}\le\sup_{|\alpha|=k}
\|\partial^{\alpha}u\|_{L^a}
\le C^{k+1}\lambda^{k}\|u\|_{L^a},\\
&&\mathrm{Supp} \,\hat{u}\subset\lambda \mathcal{C}\Rightarrow \|\sigma(D)u\|_{L^b}\le C_{\sigma,m}\lambda^{m+3(1/a-1/b)}\|u\|_{L^a}.
\end{align*}
\end{Lemma}
In order to prove the global well-poseness of our main Theorem, we also need the
following anisotropic Bernstein's Lemma .
\begin{Lemma}(see\cite{bcd})\label{bernstein2}
Let $\mathcal{B}_h$ (resp. $\mathcal{B}_v$)  be a ball of $\R_h^2$ (resp. $\R_v$) and $\mathcal{C}_h$ (resp. $\mathcal{C}_v$)a ring of $\R_h^2$ (resp. $\R_v$). Let $1\le p_2\le p_1\le\infty$ and $1\le q_2\le q_1\le\infty$. Then there hold:

If the support of $\hat{f}$ is included in $2^k\mathcal{B}_h$, then
$$\|\partial^\alpha_{x_h}f\|_{L^{p_1}_{h}(L^{q_1}_{v})}\lesssim2^{k\left(|\alpha|
+2(\frac{1}{p_2}-\frac{1}{p_1})\right)}\|f\|_{L^{p_2}_{h}(L^{q_1}_{v})}.$$

If the support of $\hat{f}$ is included in $2^l\mathcal{B}_v$, then
$$\|\partial^\beta_{x_3}f\|_{L^{p_1}_{h}(L^{q_1}_{v})}\lesssim2^{l(|\beta|
+(\frac{1}{q_2}-\frac{1}{q_1}))}\|f\|_{L^{p_1}_{h}(L^{q_2}_{v})}.$$

If the support of $\hat{f}$ is included in $2^k\mathcal{C}_h$, then
$$\|f\|_{L^{p_1}_{h}(L^{q_1}_{v})}\lesssim2^{-kN}\sup_{|\alpha|=N}\|\partial^\alpha_{x_h}f\|_{L^{p_1}_{h}(L^{q_1}_{v})}.$$

If the support of $\hat{f}$ is included in $2^l\mathcal{C}_v$, then
$$\|f\|_{L^{p_1}_{h}(L^{q_1}_{v})}\lesssim2^{-lN}\|\partial^N_3f\|_{L^{p_1}_{h}(L^{q_1}_{v})}.$$
\end{Lemma}
In the sequel, we shall frequently use Bony's decomposition from \cite{bony} in the homogeneous context:
\begin{align}\label{bony}
uv=\dot{T}_uv+\dot{T}_vu+\dot{R}(u,v)=\dot{T}_uv+\mathcal{R}(u,v),
\end{align}
where
$$\dot{T}_uv\triangleq\sum_{j\in Z}\dot{S}_{j-1}u\dot{\Delta}_jv, \hspace{0.5cm}\dot{R}(u,v)\triangleq\sum_{j\in Z}
\dot{\Delta}_ju\tilde{\Delta}_jv,$$
and
$$ \tilde{\Delta}_jv\triangleq\sum_{|j-j'|\le1}\Delta_{j'}v,\hspace{0.5cm} \mathcal{R}(u,v)\triangleq\sum_{j\in Z}\dot{S}_{j+2}v\dot{\Delta}_ju$$

In the sequel, we shall frequently use Bony's decomposition from \cite{bony} in the homogeneous context:
$$uv=\dot{T}_uv+\dot{R}(u,v)=\dot{T}_uv+\dot{T}_vu+\dot{\mathcal{R}}(u,v),$$
where
$$\dot{T}_uv\triangleq\sum_{j\in Z}\dot{S}_{j-1}u\dot{\Delta}_jv, \hspace{0.5cm}\dot{R}(u,v)\triangleq\sum_{j\in Z}
\dot{\Delta}_ju\dot{S}_{j+2}v,\hspace{0.2cm} \hspace{0.5cm}\dot{\mathcal{R}}(u,v)\triangleq\sum_{j\in Z}
\dot{\Delta}_ju\tilde{\dot{\Delta}}_jv\hspace{0.2cm} \mathrm{and}\hspace{0.2cm} \tilde{\dot{\Delta}}_jv\triangleq\sum_{|j-j'|\le1}\dot{\Delta}_{j'}v.$$
As an application of the above basic facts on Littlewood-Paley theory, we present the following
product laws in Besov spaces, which will be constantly used in the sequel.
\begin{Lemma}(see\cite{abidi2007,paicu2012})\label{daishu1}
Let $s_1\le3/2,s_2\le3/2$ with $s_1+s_2>0$. Let $a\in \dot{B}_{2,1}^{s_1}(\R^3),b\in \dot{B}_{2,1}^{s_2}(\R^3)$. Then $ab\in \dot{B}_{2,1}^{s_1+s_2-3/2}(\R^3)$, and
$$\|ab\|_{\dot{B}_{2,1}^{s_1+s_2-3/2}}\le C\|a\|_{\dot{B}_{2,1}^{s_1}}\|b\|_{\dot{B}_{2,1}^{s_2}}.$$
\end{Lemma}
\begin{Lemma}(see\cite{abidi2007})\label{daishu2}
Let $s_1\le3/2,s_2<3/2$ with $s_1+s_2\ge0$. Assume that $f\in \dot{B}_{2,1}^{s_1}(\R^3)$
 and $g\in \dot{B}_{2,\infty}^{s_2}(\R^3)$. Then there holds
$$\|fg\|_{\dot{B}_{2,\infty}^{s_1+s_2-3/2}}\le C\|f\|_{\dot{B}_{2,1}^{s_1}}\|g\|_{\dot{B}_{2,\infty}^{s_2}}.$$
\end{Lemma}

\begin{Lemma}\label{jiaohuanzi} \cite{danchincpde}
Let $s\in(-3/2,5/2]$. There exists a sequence $c_q\in l^1(\Z)$ such that $\|c_1\|_{l^1}=1$ and a constant C such that
\begin{align}
\|[v\cdot\nabla, \dot{\Delta}_q]u\|_{L^{2}}\le Cc_{q}2^{-qs}\|\nabla v\|_{\dot{B}_{2,1}^{3/2}}\|u\|_{\dot{B}_{2,1}^{s}}.
\end{align}
In the limit case $s=-3/2$, we have
$$\sup_{j\in\Z}2^{-3/2j}\|[v\cdot\nabla, \dot{\Delta}_j]u\|_{L^{2}}\le C\|\nabla v\|_{\dot{B}_{2,1}^{3/2}}\|u\|_{\dot{B}_{2,\infty}^{-3/2}}.$$
\end{Lemma}
To deal with the pressure term, we also need:
\begin{Lemma}\cite{abidi2012}\label{jiaohuanzi2}
$(\mathrm{i})$ Let $a\in B_{2,1}^{5/2}(\R^3)$ and $\nabla P\in L^2(\R^3)$ , then there holds
\begin{align}
\|[a, \Delta_q]\nabla P\|_{L^{2}}\le Cc_{q}2^{-q}\|a\|_{B_{2,1}^{5/2}}\|\nabla P\|_{L^2}.
\end{align}
$(\mathrm{ii})$ Let $a\in H^2(\R^3)$ and $\nabla P\in H^{-1}(\R^3)$, then there holds
\begin{align}
\|[a, \Delta_q]\nabla P\|_{L^{2}}\le Cc_{q}2^{q/2}\|a\|_{H^2}\|\nabla P\|_{H^{-1}}.
\end{align}
\end{Lemma}
The proof of the following two lemmas can be found in \cite{zhaicuili}, here we also  outline its proof for completeness.
\begin{Lemma}\label{guanjian}
For $2\le q,r\le\infty,v=(v^h,v^3)\in \tilde{L}^\infty((0,+\infty);\dot{B}_{2,1}^{1/2})\cap
L^1((0,+\infty);\dot{B}_{2,1}^{5/2})$ with $\mathrm{div}v=0$, there hold
 \begin{align}\label{guanjian1}
\|\dot{\Delta}_jv^3\|_{L_t^2(L^q_h(L^r_v))}\le Cd_j2^{-j(2/q+1/r)}\|v^3\|_{\tilde{L}_t^2({\dot{B}_{2,1}^{3/2}})}^{1/2+1/r}
\|v^h\|_{\tilde{L}_t^2({\dot{B}_{2,1}^{3/2}})}^{1/2-1/r},
\end{align}
and
 \begin{align}\label{guanjian2}
\|\dot{\Delta}_jv^3\|_{L_t^1(L^q_h(L^r_v))}\le Cd_j2^{-j(1+2/q+1/r)}\|v^3\|_{L_t^1({\dot{B}_{2,1}^{5/2}})}^{1/2+1/r}
\|v^h\|_{L_t^1({\dot{B}_{2,1}^{5/2}})}^{1/2-1/r}.
\end{align}
\end{Lemma}
$\mathbf{Proof:}$
Firstly, according to $\mathrm{div}v=0$ and the following version of Gagliardo-Nirenberg inequality:
 $\forall f\in\mathcal{D}(\R)$
 \begin{align}\label{}
\|f\|_{L^p(\R)}\le C \|f\|^{1/2+1/p}_{L^2(\R)}\|\nabla f\|_{L^2(\R)}^{1/2-1/p},\quad \forall\, 2\le p\le\infty,
\end{align}
we have
\begin{align}\label{}
\|\dot{\Delta}_jv^3\|_{L^r_v}&\le C\|\dot{\Delta}_jv^3\|_{L^2_v}^{1/2+1/r}\|\partial_3\dot{\Delta}_jv^3\|_{L^2_v}^{1/2-1/r}\nonumber\\
&\le C\|\dot{\Delta}_jv^3\|_{L^2_v}^{1/2+1/r}\|\dot{\Delta}_j(\mathrm{div}_hv^h)\|_{L^2_v}^{1/2-1/r}
\le C2^{j(1/2-1/r)}\|\dot{\Delta}_jv^3\|_{L^2_v}^{1/2+1/r}\|\dot{\Delta}_jv^h\|_{L^2_v}^{1/2-1/r}.
\end{align}
Thus,
by Lemma \ref{bernstein}, we have
\begin{align}\label{}
\|\dot{\Delta}_jv^3\|_{L_t^2(L^q_h(L^r_v))}&\le C2^{-j(2/q-1)}\|\dot{\Delta}_jv^3\|_{L_t^2(L^2_h(L^r_v))}\nonumber\\
&\le C2^{j(3/2-2/q-1/r)}\|\dot{\Delta}_jv^3\|_{L_t^2(L^2_h(L^2_v))}^{1/2+1/r}
\|\dot{\Delta}_jv^h\|_{L_t^2(L^2_h(L^2_v))}^{1/2-1/r}\nonumber\\
&\le Cd_j2^{-j(2/q+1/r)}\|v^3\|_{\tilde{L}_t^2({\dot{B}_{2,1}^{3/2}})}^{1/2+1/r}
\|v^h\|_{\tilde{L}_t^2({\dot{B}_{2,1}^{3/2}})}^{1/2-1/r}.
\end{align}
Similarly, we have (\ref{guanjian2}).
From the above Lemmas, we can obtain the following important Lemma.
\begin{Lemma}\label{zhongyaoguji}
Let $u\in\tilde{L}_t^\infty({\dot{B}_{2,1}^{1/2}})\cap L_t^1({\dot{B}_{2,1}^{5/2}})$ and $\mathrm{div}\,u=0$, one has
\begin{align}
&2^j\|\dot{\Delta}_j(u^3u^h)\|_{L_t^1(L^2)}+\|\dot{\Delta}_j(u^3\mathrm{div}_hu^h)\|_{L_t^1(L^2)}\nonumber\\
\le& Cd_j2^{-j/2}\Bigg(
\|u^3\|_{{L}_{t}^1({\dot{B}_{2,1}^{5/2}})}^{1/2}
\|u^h\|_{{L}_{t}^1({\dot{B}_{2,1}^{5/2}})}^{1/2}\|u^h\|_{\tilde{L}_t^\infty({\dot{B}_{2,1}^{1/2}})}
\nonumber\\
&+(\|u^3\|_{{L}_{t}^1({\dot{B}_{2,1}^{5/2}})}+\|u^3\|_{\tilde{L}_t^\infty({\dot{B}_{2,1}^{1/2}})})^{1/2}
(\|u^h\|_{{L}_{t}^1({\dot{B}_{2,1}^{5/2}})}+\|u^h\|_{\tilde{L}_t^\infty({\dot{B}_{2,1}^{1/2}})})^{3/2}
\Bigg)\nonumber\\
\triangleq& Cd_j2^{-j/2}F(u^3,u^h).
\end{align}
\end{Lemma}
$\mathbf{Proof:}$ Firstly,
thanks to Bony's decomposition, we have
\begin{align}\label{}
u^3u^h+u^3\mathrm{div}_hu^h
=&\dot{T}_{u^3}u^h+\dot{T}_{u^h}u^3+\dot{T}_{u^3}\mathrm{div}_hu^h+\dot{T}_{\mathrm{div}_hu^h}u^3\nonumber\\
&+\dot{R}(u^3,u^h)+\dot{R}(u^3,\mathrm{div}_hu^h).
\end{align}
Applying Lemma \ref{bernstein} and (\ref{guanjian1}) with $q=r=\infty$ gives rise to
\begin{align}\label{Q1}
&2^j\|\dot{\Delta}_j(\dot{T}_{u^3}u^h)\|_{L_t^1(L^2)}+\|\dot{\Delta}_j(\dot{T}_{u^3}\mathrm{div}_hu^h)\|_{L_t^1(L^2)}\nonumber\\
&\le C\sum_{|j'-j|\le5}\left(2^j\|\dot{S}_{j'-1}u^3\|_{L_t^2(L^\infty)}\|\dot{\Delta}_{j'}u^h\|_{L_t^2(L^2)}+ \|\dot{S}_{j'-1}u^3\|_{L_t^2(L^\infty)}\|\dot{\Delta}_{j'}(\mathrm{div}_hu^h)\|_{L_t^2(L^2)}\right)\nonumber\\
&\le C\sum_{|j'-j|\le5}\sum_{j''\le j'-2}\left(2^j\|\dot{\Delta}_{j''}u^3\|_{L_t^2(L^\infty)}\|\dot{\Delta}_{j'}u^h\|_{L_t^2(L^2)}+ \|\dot{\Delta}_{j''}u^3\|_{L_t^2(L^\infty)}\|\dot{\Delta}_{j'}(\mathrm{div}_hu^h)\|_{L_t^2(L^2)}\right)\nonumber\\
&\le Cd_j 2^{-j/2}\|u^3\|_{\tilde{L}_t^2({\dot{B}_{2,1}^{{3}/{2}}})}^{1/2}\|u^h\|^{3/2}_{{\tilde{L}}_{t}^2({\dot{B}_{2,1}^{{3}/{2}}})}
\nonumber\\
&\le Cd_j 2^{-j/2}(\|u^3\|_{{L}_{t}^1({\dot{B}_{2,1}^{5/2}})}+\|u^3\|_{\tilde{L}_t^\infty({\dot{B}_{2,1}^{1/2}})})^{{1}/{2}}
(\|u^h\|_{{L}_{t}^1({\dot{B}_{2,1}^{5/2}})}+\|u^h\|_{\tilde{L}_t^\infty({\dot{B}_{2,1}^{1/2}})})^{{3}/{2}}.
\end{align}
Similarly, we have
\begin{align}\label{Q2}
&\|\dot{\Delta}_j(\dot{T}_{\mathrm{div}_hu^h}u^3)\|_{L_t^1(L^2)}\nonumber\\
&\le C\sum_{|j'-j|\le5}\left(\|\dot{S}_{j'-1}(\mathrm{div}_hu^h)\|_{L_t^2(L_h^\infty(L^2_v))}\|\dot{\Delta}_{j'}u^3\|_{L_t^2(L^2_h(L^\infty_v))}\right)\nonumber\\
&\le C\sum_{|j'-j|\le5}\sum_{j''\le j'-2}\left(2^{j''} \|\dot{\Delta}_{j''}(\mathrm{div}_hu^h)\|_{L_t^2(L^2)}\|\dot{\Delta}_{j'}u^3\|_{L_t^2(L^2_h(L^\infty_v))}\right)\nonumber\\
&\le Cd_j 2^{-j/2}\|u^3\|_{\tilde{L}_t^2({\dot{B}_{2,1}^{{3}/{2}}})}^{{1}/{2}}\|u^h\|^{{3}/{2}}_{{\tilde{L}}_{t}^2({\dot{B}_{2,1}^{{3}/{2}}})}\nonumber\\
&\le Cd_j 2^{-j/2}(\|u^3\|_{{L}_{t}^1({\dot{B}_{2,1}^{5/2}})}+\|u^3\|_{\tilde{L}_t^\infty({\dot{B}_{2,1}^{1/2}})})^{{1}/{2}}
(\|u^h\|_{{L}_{t}^1({\dot{B}_{2,1}^{5/2}})}+\|u^h\|_{\tilde{L}_t^\infty({\dot{B}_{2,1}^{1/2}})})^{{3}/{2}}.
\end{align}
And while (\ref{guanjian2}) applied with $q=2,r=\infty$ gives
\begin{align}\label{Q3}
2^j\|\dot{\Delta}_j(\dot{T}_{u^h}u^3)\|_{L_t^1(L^2)} \le& C\sum_{|j'-j|\le5}\left(2^j\|\dot{S}_{j'-1}u^h\|_{L_t^\infty(L_h^\infty(L^2_v))}\|\dot{\Delta}_{j'}u^3\|_{L_t^1(L_h^2(L_v^{\infty}))}\right)\nonumber\\
\le& C\sum_{|j'-j|\le5}\sum_{j''\le j'-2}\left(2^{j+j''} \|\dot{\Delta}_{j''}u^h\|_{L_t^\infty(L^2)}\|\dot{\Delta}_{j'}u^3\|_{L_t^1(L_h^2(L_v^{\infty}))}\right)\nonumber\\
\le &Cd_j2^{-j/2}
\|u^3\|_{{L}_{t}^1({\dot{B}_{2,1}^{5/2}})}^{1/2}\|u^h\|_{{L}_{t}^1({\dot{B}_{2,1}^{5/2}})}^{1/2}
\|u^h\|_{\tilde{L}_t^\infty({\dot{B}_{2,1}^{1/2}})}
\end{align}
and
\begin{align}\label{Q4}
2^j\|\dot{\Delta}_j\dot{R}(u^3,u^h)\|_{L_t^1(L^2)} \le& C\sum_{j'\ge j-N_0}\left(2^j\|\tilde{\dot{\Delta}}_{j'}u^h\|_{L_t^\infty(L_h^\infty(L^2_v))}\|\dot{\Delta}_{j'}u^3\|_{L_t^1(L_h^2(L_v^{\infty}))}\right)\nonumber\\
\le& C\sum_{j'\ge j-N_0}\left(2^{j+j'} \|\tilde{\dot{\Delta}}_{j'}u^h\|_{L_t^\infty(L^2)}\|\dot{\Delta}_{j'}u^3\|_{L_t^1(L_h^2(L_v^{\infty}))}\right)\nonumber\\
\le &Cd_j2^{-j/2}
\|u^3\|_{{L}_{t}^1({\dot{B}_{2,1}^{5/2}})}^{{1}/{2}}\|u^h\|_{{L}_{t}^1({\dot{B}_{2,1}^{5/2}})}^{{1}/{2}}
\|u^h\|_{\tilde{L}_t^\infty({\dot{B}_{2,1}^{1/2}})}.
\end{align}
In what follows, we give the estimate of the remaining term $\dot{R}(u^3,\mathrm{div}_hu^h)$. Applying (\ref{guanjian1}) with $q=2,r=\infty$ gives
\begin{align}\label{Q5}
&\|\dot{\Delta}_j\dot{R}(u^3,\mathrm{div}_hu^h)\|_{L_t^1(L^2)}\nonumber\\
&\le C\sum_{j'\ge j-N_0}\left(2^{j}\|\tilde{\dot{\Delta}}_{j'}(\mathrm{div}_hu^h)\dot{\Delta}_{j'}u^3\|_{L_t^1(L_h^{1}(L_v^2))}\right)\nonumber\\
&\le C\sum_{j'\ge j-N_0}\left(2^{j}\|\tilde{\dot{\Delta}}_{j'}(\mathrm{div}_hu^h)\|_{L_t^2(L^2)}\|\dot{\Delta}_{j'}u^3\|_{L_t^2(L_h^2(L_v^{\infty}))}\right)\nonumber\\
&\le Cd_j 2^{-j/2}\|u^3\|_{\tilde{L}_t^2({\dot{B}_{2,1}^{{3}/{2}}})}^{1/2}\|u^h\|^{3/2}_{{\tilde{L}}_{t}^2({\dot{B}_{2,1}^{{3}/{2}}})}
\nonumber\\
&\le Cd_j 2^{-j/2}(\|u^3\|_{{L}_{t}^1({\dot{B}_{2,1}^{5/2}})}+\|u^3\|_{\tilde{L}_t^\infty({\dot{B}_{2,1}^{1/2}})})^{{1}/{2}}
(\|u^h\|_{{L}_{t}^1({\dot{B}_{2,1}^{5/2}})}+\|u^h\|_{\tilde{L}_t^\infty({\dot{B}_{2,1}^{1/2}})})^{{3}/{2}}.
\end{align}

Combining with the estimates (\ref{Q2})-(\ref{Q5}), we can finally get
\begin{align}\label{1p}
&2^j\|\dot{\Delta}_j(u^3u^h)\|_{L_t^1(L^2)}+\|\dot{\Delta}_j(u^3\mathrm{div}_hu^h)\|_{L_t^1(L^2)}\nonumber\\
\le& Cd_j2^{{-j/2}}\Bigg(
\|u^3\|_{{L}_{t}^1({\dot{B}_{2,1}^{5/2}})}^{1/2}
\|u^h\|_{{L}_{t}^1({\dot{B}_{2,1}^{5/2}})}^{1/2}\|u^h\|_{\tilde{L}_t^\infty({\dot{B}_{2,1}^{1/2}})}
\nonumber\\
&+(\|u^3\|_{{L}_{t}^1({\dot{B}_{2,1}^{5/2}})}+\|u^3\|_{\tilde{L}_t^\infty({\dot{B}_{2,1}^{1/2}})})^{{1}/{2}}
(\|u^h\|_{{L}_{t}^1({\dot{B}_{2,1}^{5/2}})}+\|u^h\|_{\tilde{L}_t^\infty({\dot{B}_{2,1}^{1/2}})})^{{3}/{2}}
\Bigg).
\end{align}
The proof of Lemma \ref{guanjian} is complete.\quad$\square$

\section{Estimates of the transport equation }
Let us first recall standard estimates in Besov spaces for the following linear transport equation:

\begin{equation}\label{3.0}
\left\{
 \begin{array}{ll}
 \partial_t f +  v \cdot\nabla f =g, \\
f(x,0)= f_0.
 \end{array} \right.
  \end{equation}
\begin{proposition} (see \cite{abidi2012})\label{zhiliangguji0}
 Let v be a divergence free vector field with $\nabla v\in L^1([0,T];\dot{B}_{2,1}^{3/2}(\R^3))$. For $s\in (-5/2, 5/2]$, given $f_0\in B_{2,1}^{s}(\R^3), g\in L^1([0,T];B_{2,1}^{s}(\R^3))$, then the system (\ref{3.0}) has a unique solution $f\in C([0,T];B_{2,1}^{s}(\R^3))$ . Moreover, there exists a constant $C$ such that for $\forall t\in[0,T]$
 \begin{align}
\|f\|_{\widetilde{L}_t^{\infty}(B_{2,1}^{s})}&\le \|f_0\|_{B_{2,1}^{s}}+C\int_0^t\|f\|_{B_{2,1}^{s}}\|\nabla v\|_{\dot{B}_{2,1}^{3/2}}d\tau+C\|g\|_{L_t^1({B_{2,1}^{s}})}.
\end{align}
\end{proposition}

\begin{proposition}\cite{abidi2012}  \label{shuyun}
Let $m\in\Z, a_0\in B_{2,1}^{3/2}(\R^3),\nabla u\in L_T^1({\dot{B}_{2,1}^{3/2}})$ with $\mathrm{div} u=0$, and $a\in C([0, T ]; B_{2,1}^{3/2}(\R^3)$ such that $(a,u)$ solves
\begin{equation}\label{3.1}
   \left\{
 \begin{array}{ll}
 \partial_t a +  u \cdot\nabla a =0, \\
a(x,0)= a_0.
 \end{array} \right.
  \end{equation}
Then there hold for $\forall t\le T$
\begin{align}\label{3.2}
\|a\|_{\tilde{L}^\infty_t(B_{2,1}^{3/2})}\le\|a_0\|_{B_{2,1}^{3/2}}e^{CU(t)},
\end{align}
\begin{align}\label{3.3}
\|a-S_ma\|_{\tilde{L}^\infty_t(B_{2,1}^{3/2})}\le\sum_{q\ge m} 2^{3q/2}\|\Delta_qa_0\|_{L^2}+\|a_0\|_{B_{2,1}^{3/2}}(e^{CU(t)}-1),
\end{align}
with $U(t)=\|\nabla u\|_{L^1_t(\dot{B}_{2,1}^{3/2})}$.
\end{proposition}
Finally, we give estimates for the following constant coefficient parabolic system:
\begin{equation}\label{refangcheng}
   \left\{
 \begin{array}{ll}
 \partial_t u-\sigma\Delta u  =f, \\
u(x,0)= u_0.
 \end{array} \right.
  \end{equation}
\begin{proposition} (see \cite{danchincpde}£© \label{ziyoufangcheng}
Assume that $\sigma\ge0$, then there exists a universal constant $\kappa$ such that for all $s\in\R$ and $T\in\R^+,$
$$\|u\|_{\tilde{L}^\infty_T(\dot{B}_{2,1}^s)}\le \|u_0\|_{\dot{B}_{2,1}^s}+\|f\|_{L_T^1({\dot{B}_{2,1}^{s}})},
$$
$$\kappa\|u\|_{L_T^1({\dot{B}_{2,1}^{s+2}})}\le\sum_{q\in\Z}2^{qs}(1-e^{-\kappa\sigma2^{2q}T})
(\|\dot{\Delta}_qu_0\|_{L^2}+\|\dot{\Delta}_qf\|_{L_T^1(L^2)}).
$$
\end{proposition}

\begin{center}
\section{Local wellposedness of Theorem 1.1}
\end{center}

\subsection{Existence of Theorem 1.1}
The proof of existence of a solution is performed in a standard manner. We begin
by solving an approximate problem and we prove that the solutions are uniformly
bounded. The last step consists in studying the convergence to a solution of the
initial equation by a compactness argument.

Step 1: Construction of smooth approximate solutions.

We first smooth out the initial data. For $n \in\N$, let
$a_0^n=\dot{S}_na_0-\dot{S}_{-n}a_0,\,u_0^n=\dot{S}_nu_0-\dot{S}_{-n}u_0$, $d_0^n=\dot{S}_nd_0-\dot{S}_{-n}d_0$  then ($a_0^n,u_0^n,d_0^n)\in H^\infty(\R^3)$. Now, owing to \cite{wenhuanyao}, we deduce that system (\ref{1.1}) with the initial data
$(a_0^n,u_0^n,d_0^n)$ admits a unique local in time solution $(a^n, u^n,d^n,\nabla  P^n)$ verifying
$$a^n\in C([0,T^n];H^{s+1}(\R^3)),$$
$$u^n\in C([0,T^n];H^{s+1}(\R^3))\cap L^1_{loc}([0,T^n];H^{s+2}(\R^3)),$$
$$d^n\in C([0,T^n];H^{s+2}(\R^3))\cap L^1_{loc}([0,T^n];H^{s+3}(\R^3))$$
and
$$\nabla P^n\in L^1([0,T^n];H^{s}(\R^3)) \,\,{\rm with}\,\, s>1/2.$$

Step 2: Uniform estimates to the approximate solutions.

Our goal is to prove that there exists a positive time $0 < T < \inf T^n(n\in\N),$
such that $(a^n, u^n,d^n,\nabla P^n)$ is uniformly bounded in the space
\begin{align*}
E_T=&\tilde{L}_T^\infty({\dot{B}_{2,1}^{3/2}})\times\tilde{L}_T^\infty({\dot{B}_{2,1}^{1/2}})\cap L_T^1({\dot{B}_{2,1}^{5/2}})\times\tilde{L}_T^\infty({\dot{B}_{2,1}^{3/2}})\cap L_T^1({\dot{B}_{2,1}^{7/2}})\times L_T^1({\dot{B}_{2,1}^{1/2}}).
\end{align*}
Let $u^n=u_L^n+\bar{u}^n,d^n=d_L^n+\bar{d}^n$ which $u_L^n=e^{t\Delta }u_0^n,d_L^n=e^{t\Delta }d_0^n$, then
$(a^n, \bar{u}^n,\bar{d}^n,\nabla P^n)$ solves
\begin{eqnarray}\label{feixianxing}
\left\{\begin{aligned}
&\partial_ta^n+(u_L^n+\bar{u}^n)\cdot\nabla a^n=0,\\
&\partial_t\bar{u}^n+ u_L^n\cdot \nabla \bar{u}^n+(1+a^n)(\nabla P^n-\Delta \bar{u}^n)=M_n,\\
&\partial_t\bar{d}^n-\Delta \bar{d}^n+u_L^n \cdot \nabla \bar{d}^n=G_n,\\
&\mathrm{\mathrm{div}}\bar{u}^n=0,\\
&(a^n,\bar{u}^n,\bar{d}^n)|_{t=0}=(a_0,0,0),
\end{aligned}\right.
\end{eqnarray}
with
\begin{align*}
M_n=&-u_L^n\cdot \nabla u_L^n-\bar{u}^n\cdot \nabla u_L^n-\bar{u}^n\cdot \nabla \bar{u}^n+a^n\Delta u_L^n
\\&+(1+a^n)(\nabla d_L^n\cdot \Delta d_L^n+\nabla d_L^n\cdot \Delta \bar{d}^n+\nabla \bar{d}^n\cdot \Delta \bar{d}^n+\nabla \bar{d}^n\cdot \Delta d_L^n)
\end{align*}
and
\begin{align*}
G_n=&-u_L^n\cdot \nabla d_L^n-\bar{u}^n\cdot \nabla d_L^n-\bar{u}^n\cdot \nabla \bar{d}^n+|\nabla d_L^n|^2d_L^n\\
&+|\nabla d_L^n|^2\bar{d}^n+|\nabla \bar{d}^n|^2d_L^n+|\nabla \bar{d}^n|^2\bar{d}^n
+2\nabla d_L^n\cdot\nabla \bar{d}^nd_L^n+2\nabla d_L^n\cdot\nabla \bar{d}^n\bar{d}^n.
\end{align*}

It is easy  to observe from Proposition \ref{ziyoufangcheng} that
\begin{align}\label{4.0}
\|u_L^n\|_{\tilde{L}^\infty_t(\R^+;\dot{B}_{2,1}^{1/2})}\le C\|u_0\|_{\dot{B}_{2,1}^{1/2}},\hspace{0.5cm}
\|d_L^n\|_{\tilde{L}^\infty_t(\R^+;\dot{B}_{2,1}^{3/2})}\le C\|d_0\|_{\dot{B}_{2,1}^{3/2}}
\end{align}
and
\begin{align}\label{4.1}
&\|u_L^n\|_{L_t^1({\dot{B}_{2,1}^{5/2}})}\le C\sum_{q\in\Z}2^{q/2}(1-e^{-ct2^{2q}})
\|\dot{\Delta}_qu_0\|_{L^2},\nonumber\\
&\|d_L^n\|_{L_t^1({\dot{B}_{2,1}^{7/2}})}\le C\sum_{q\in\Z}2^{{3q}/2}(1-e^{-ct2^{2q}})
\|\dot{\Delta}_qd_0\|_{L^2}.
\end{align}
Let
\begin{align*}
Z^n(t)\triangleq&\|\bar{u}^n\|_{\tilde{L}^\infty_t(\dot{B}_{2,1}^{1/2})}+\|\bar{d}^n\|_{\tilde{L}^\infty_t(\dot{B}_{2,1}^{3/2})}
+\|\bar{u}^n\|_{L_t^1({\dot{B}_{2,1}^{5/2}})}+\|\bar{d}^n\|_{L_t^1({\dot{B}_{2,1}^{7/2}})}
+\|\nabla P^n\|_{L_t^1({\dot{B}_{2,1}^{1/2}})}.
\end{align*}
Firstly, we give the estimate of d.

Applying $\dot{\Delta}_q$ to the third equation of (\ref{feixianxing}), and then a standard commutator's process, gives
\begin{align}\label{4.2}
&\partial_t\dot{\Delta}_q\bar{d}^n+u_L^n \cdot \nabla\dot{\Delta}_q\bar{d}^n+\bar{u}^n \cdot \nabla\dot{\Delta}_q\bar{d}^n-\Delta \dot{\Delta}_q\bar{d}^n \nonumber\\
=&
[u_L^n \cdot \nabla,\dot{\Delta}_q]\bar{d}^n+[\bar{u}^n \cdot \nabla,\dot{\Delta}_q]\bar{d}^n+\dot{\Delta}_q(|\nabla d_L^n|^2d_L^n)\nonumber\\
&+\dot{\Delta}_q(|\nabla d_L^n|^2\bar{d}^n)+\dot{\Delta}_q(|\nabla \bar{d}^n|^2d_L^n)
+\dot{\Delta}_q(|\nabla \bar{d}^n|^2\bar{d}^n)
\nonumber\\
&+\dot{\Delta}_q(2\nabla d_L^n\cdot\nabla \bar{d}^nd_L^n)
+\dot{\Delta}_q(2\nabla d_L^n\cdot\nabla \bar{d}^n\bar{d}^n)+\dot{\Delta}_q(u_L^n \cdot \nabla d_L^n)+\dot{\Delta}_q(\bar{u}^n \cdot \nabla d_L^n).
\end{align}
Taking the $L^2$ inner product of the above resulting
equation with $\dot{\Delta}_q\bar{d}^n$ and integrating with respect to $t$, we can finally get
\begin{align}\label{4.3}
&\|\bar{d}^n\|_{\tilde{L}^\infty_t(\dot{B}_{2,1}^{3/2})}+
\|\bar{d}^n\|_{L^1_t(\dot{B}_{2,1}^{7/2})}\nonumber\\
\lesssim&\sum_{q\in\Z}2^{{3q}/2} \|[u_L^n \cdot \nabla,\dot{\Delta}_q]\bar{d}^n\|_{L^1_t(L^2)}+\sum_{q\in\Z}2^{{3q}/2} \|[\bar{u}^n \cdot \nabla,\dot{\Delta}_q]\bar{d}^n\|_{L^1_t(L^2)}\nonumber\\
&+\|u_L^n \cdot \nabla d_L^n\|_{L_t^1({\dot{B}_{2,1}^{3/2}})}
+\|\bar{u}^n \cdot \nabla d_L^n\|_{L_t^1({\dot{B}_{2,1}^{3/2}})}
+\||\nabla d_L^n|^2d_L^n\|_{L_t^1({\dot{B}_{2,1}^{3/2}})}
+\||\nabla d_L^n|^2\bar{d}^n\|_{L_t^1({\dot{B}_{2,1}^{3/2}})}\nonumber\\
&+\||\nabla \bar{d}^n|^2d_L^n\|_{L_t^1({\dot{B}_{2,1}^{3/2}})}
+\||\nabla \bar{d}^n|^2\bar{d}^n\|_{L_t^1({\dot{B}_{2,1}^{3/2}})}
+\|\nabla d_L^n\cdot\nabla \bar{d}^nd_L^n\|_{L_t^1({\dot{B}_{2,1}^{3/2}})}
+\|\nabla d_L^n\cdot\nabla \bar{d}^n\bar{d}^n\|_{L_t^1({\dot{B}_{2,1}^{3/2}})}.
\end{align}
Applying Lemma \ref{daishu1}, \ref{jiaohuanzi} yields that
\begin{align}\label{d1}
\sum_{q\in\Z}2^{{3q}/2} \|[u_L^n \cdot \nabla,\dot{\Delta}_q]\bar{d}^n\|_{L^1_t(L^2)}\le C \int_0^t\|u_L^n\|_{\dot{B}_{2,1}^{5/2}}\|\bar{d}^n\|_{\dot{B}_{2,1}^{3/2}}d\tau\le C\int_0^t\|u_L^n\|_{\dot{B}_{2,1}^{5/2}}Z^n(\tau)d\tau,
\end{align}
\begin{align}\label{d2}
\sum_{q\in\Z}2^{{3q}/2} \|[\bar{u}^n \cdot \nabla,\dot{\Delta}_q]\bar{d}^n\|_{L^1_t(L^2)}\le C
\int_0^t\|\bar{u}^n\|_{\dot{B}_{2,1}^{5/2}}\|\bar{d}^n\|_{\dot{B}_{2,1}^{3/2}}d\tau\le C(Z^n(t))^2,
\end{align}
\begin{align}\label{d3}
\|u_L^n \cdot \nabla d_L^n\|_{L_t^1({\dot{B}_{2,1}^{3/2}})}\le & C\int_0^t\|u_L^n\|_{\dot{B}_{2,1}^{3/2}}\|d_L^n\|_{\dot{B}_{2,1}^{5/2}}d\tau\le C
\int_0^t\|u_L^n\|_{\dot{B}_{2,1}^{1/2}}^{1/2}\|u_L^n\|_{\dot{B}_{2,1}^{5/2}}^{1/2}
\|d_L^n\|_{\dot{B}_{2,1}^{3/2}}^{1/2}\|d_L^n\|_{\dot{B}_{2,1}^{7/2}}^{1/2}d\tau\nonumber\\
\le &C\int_0^t\|u_L^n\|_{\dot{B}_{2,1}^{1/2}}\|u_L^n\|_{\dot{B}_{2,1}^{5/2}}d\tau+
C\int_0^t\|d_L^n\|_{\dot{B}_{2,1}^{3/2}}\|d_L^n\|_{\dot{B}_{2,1}^{7/2}}d\tau\nonumber\\
\le &C (\|u_0\|_{\dot{B}_{2,1}^{1/2}}+\|d_0\|_{\dot{B}_{2,1}^{3/2}} )( \|u_L^n\|_{L^1_t(\dot{B}_{2,1}^{5/2})} +\|d_L^n\|_{L^1_t(\dot{B}_{2,1}^{7/2})} ),
\end{align}
\begin{align}\label{d4}
\|\bar{u}^n \cdot \nabla d_L^n\|_{L_t^1({\dot{B}_{2,1}^{3/2}})}\le &
C\int_0^t\|\bar{u}^n\|_{\dot{B}_{2,1}^{3/2}}\|d_L^n\|_{\dot{B}_{2,1}^{5/2}} d\tau\le C
\int_0^t\|\bar{u}^n\|_{\dot{B}_{2,1}^{1/2}}^{1/2}\|\bar{u}^n\|_{\dot{B}_{2,1}^{5/2}}^{1/2}
\|d_L^n\|_{\dot{B}_{2,1}^{3/2}}^{1/2}\|d_L^n\|_{\dot{B}_{2,1}^{7/2}}^{1/2}d\tau\nonumber\\
\le &C (Z^n(t))^2+C \|d_0\|_{\dot{B}_{2,1}^{3/2}}\|d_L^n\|_{L^1_t(\dot{B}_{2,1}^{7/2})},
\end{align}
\begin{align}\label{d5}
\||\nabla d_L^n|^2d_L^n\|_{L_t^1({\dot{B}_{2,1}^{3/2}})}\le  C\int_0^t\|d_L^n\|_{\dot{B}_{2,1}^{3/2}}\|d_L^n\|_{\dot{B}_{2,1}^{5/2}}^2d\tau\le  C\int_0^t\|d_L^n\|_{\dot{B}_{2,1}^{3/2}}^2\|d_L^n\|_{\dot{B}_{2,1}^{7/2}}d\tau\le C \|d_0\|_{\dot{B}_{2,1}^{3/2}}^2\|d_L^n\|_{L^1_t(\dot{B}_{2,1}^{7/2})},
\end{align}
\begin{align}\label{d6}
\||\nabla d_L^n|^2\bar{d}^n\|_{L_t^1({\dot{B}_{2,1}^{3/2}})}\le &C\int_0^t\|\bar{d}^n\|_{\dot{B}_{2,1}^{3/2}}\|d_L^n\|_{\dot{B}_{2,1}^{5/2}}^2d\tau\nonumber\\
\le &C\int_0^t\|\bar{d}^n\|_{\dot{B}_{2,1}^{3/2}}\|d_L^n\|_{\dot{B}_{2,1}^{3/2}}\|d_L^n\|_{\dot{B}_{2,1}^{7/2}}d\tau
\le
C\int_0^t\|d_0\|_{\dot{B}_{2,1}^{3/2}}\|d_L^n\|_{\dot{B}_{2,1}^{7/2}}Z^n(\tau)d\tau,
\end{align}
\begin{align}\label{d7}
\||\nabla \bar{d}^n|^2d_L^n\|_{L_t^1({\dot{B}_{2,1}^{3/2}})}
\le& C\int_0^t\|d_L^n\|_{\dot{B}_{2,1}^{3/2}}\|\bar{d}^n\|_{\dot{B}_{2,1}^{5/2}}^2d\tau\nonumber\\
\le& C\int_0^t\|d_L^n\|_{\dot{B}_{2,1}^{3/2}}\|\bar{d}^n\|_{\dot{B}_{2,1}^{3/2}}\|\bar{d}^n\|_{\dot{B}_{2,1}^{7/2}}d\tau
\le C \|d_0\|_{\dot{B}_{2,1}^{3/2}} (Z^n(t))^2,
\end{align}
\begin{align}\label{d8}
\||\nabla \bar{d}^n|^2\bar{d}^n\|_{L_t^1({\dot{B}_{2,1}^{3/2}})}
\le& C\int_0^t\|\bar{d}^n\|_{\dot{B}_{2,1}^{3/2}}\|\bar{d}^n\|_{\dot{B}_{2,1}^{5/2}}^2d\tau\nonumber\\
\le& C\int_0^t\|\bar{d}^n\|_{\dot{B}_{2,1}^{3/2}}^2\|\bar{d}^n\|_{\dot{B}_{2,1}^{7/2}}d\tau
\le C (Z^n(t))^3,
\end{align}

\begin{align}\label{d9}
\|\nabla d_L^n\cdot\nabla \bar{d}^nd_L^n\|_{L_t^1({\dot{B}_{2,1}^{3/2}})}
\le& C\int_0^t  \|d_L^n\|_{\dot{B}_{2,1}^{5/2}} \|\bar{d}^n\|_{\dot{B}_{2,1}^{5/2}} \|d_L^n\|_{\dot{B}_{2,1}^{3/2}}d\tau\nonumber\\
\le& C\int_0^t  \|\bar{d}^n\|_{\dot{B}_{2,1}^{3/2}}^{1/2}\|\bar{d}^n\|_{\dot{B}_{2,1}^{7/2}}^{1/2} \|d_L^n\|_{\dot{B}_{2,1}^{3/2}}^{3/2}\|d_L^n\|_{\dot{B}_{2,1}^{7/2}}^{1/2} d\tau\nonumber\\
\le& C\int_0^t  \|\bar{d}^n\|_{\dot{B}_{2,1}^{3/2}}\|\bar{d}^n\|_{\dot{B}_{2,1}^{7/2}}d\tau+ C\int_0^t \|d_L^n\|_{\dot{B}_{2,1}^{3/2}}^{3}\|d_L^n\|_{\dot{B}_{2,1}^{7/2}} d\tau\nonumber\\
\le& C (Z^n(t))^2+ C\|d_0\|_{\dot{B}_{2,1}^{3/2}}^3\|d_L^n\|_{L^1_t(\dot{B}_{2,1}^{7/2})},
\end{align}
\begin{align}\label{d10}
\|\nabla d_L^n\cdot\nabla \bar{d}^n\bar{d}^n\|_{L_t^1({\dot{B}_{2,1}^{3/2}})}
\le& C\int_0^t  \|d_L^n\|_{\dot{B}_{2,1}^{5/2}} \|\bar{d}^n\|_{\dot{B}_{2,1}^{5/2}} \|\bar{d}^n\|_{\dot{B}_{2,1}^{3/2}}d\tau\nonumber\\
\le& C\int_0^t  \|d_L^n\|_{\dot{B}_{2,1}^{3/2}}^{1/2}\|d_L^n\|_{\dot{B}_{2,1}^{7/2}}^{1/2}
\|\bar{d}^n\|_{\dot{B}_{2,1}^{3/2}}^{3/2}\|\bar{d}^n\|_{\dot{B}_{2,1}^{7/2}}^{1/2}d\tau\nonumber\\
\le& C\int_0^t  \|d_L^n\|_{\dot{B}_{2,1}^{3/2}}\|d_L^n\|_{\dot{B}_{2,1}^{7/2}}d\tau
+C\int_0^t\|\bar{d}^n\|_{\dot{B}_{2,1}^{3/2}}^{3}\|\bar{d}^n\|_{\dot{B}_{2,1}^{7/2}}d\tau\nonumber\\
\le& C(Z^n(t))^4+\|d_0\|_{\dot{B}_{2,1}^{3/2}}\|d_L^n\|_{L^1_t(\dot{B}_{2,1}^{7/2})}.
\end{align}
Taking above estimates into (\ref{4.3}), we have
\begin{align}\label{d11}
&\|\bar{d}^n\|_{\tilde{L}^\infty_t(\dot{B}_{2,1}^{3/2})}+\|\bar{d}^n\|_{L^1_t(\dot{B}_{2,1}^{7/2})}\nonumber\\
\le& C (1+\|d_0\|_{\dot{B}_{2,1}^{3/2}})(Z^n(t))^2+C(Z^n(t))^4
+C\int_0^t(\|u_L^n\|_{\dot{B}_{2,1}^{5/2}} +(1+\|d_0\|_{\dot{B}_{2,1}^{3/2}})\|d_L^n\|_{\dot{B}_{2,1}^{7/2}})Z^n(\tau)d\tau\nonumber\\
&+C \left(\|u_0\|_{\dot{B}_{2,1}^{1/2}}+(1+\|d_0\|_{\dot{B}_{2,1}^{3/2}}^2)(1+\|d_0\|_{\dot{B}_{2,1}^{3/2}})\right)( \|u_L^n\|_{L^1_t(\dot{B}_{2,1}^{5/2})} +\|d_L^n\|_{L^1_t(\dot{B}_{2,1}^{7/2})} ).
\end{align}

In what follows, we mainly discuss the momentum equation. In the general case of large perturbations of a constant density state, solving the nonconstant coefficients elliptic equation may be a complicated problem.
Here, we follow the Proposition 6 in  \cite{abidi2012}   to give the following key proposition:
\begin{proposition}\label{dongliang}
Let $a\in L^\infty_T(B_{2,1}^{3/2})(\R^3),u_0\in \dot{B}_{2,1}^{3/2}(\R^3)$ which satisfy
$1+a\ge \underline{b}$ for some positive constant $\underline{b}$ and $\mathrm{div}u_0=0$. Let $f$ belong to $L_T^1({\dot{B}_{2,1}^{1/2}})$ with $\mathrm{div} f\in L^1_T(\dot{H}^{-1}),v\in L^\infty_T(\dot{B}_{2,1}^{1/2}),\nabla v\in L_T^1({\dot{B}_{2,1}^{3/2}})$ and $(u,\nabla P)\in C([0, T ]; \dot{B}_{2,1}^{1/2}\cap L^1_{loc}((0,T);\dot{B}_{2,1}^{5/2})\times L^1_{loc}((0,T);\dot{B}_{2,1}^{1/2})$, which solve
\begin{equation}\label{4.300}
\left\{
 \begin{array}{ll}
 \partial_tu+v \cdot \nabla u-(1+a)(\Delta u-\nabla P)=f, \\
 \mathrm{div} u=0,  \\
u(x,0)=u_0.
 \end{array} \right.
 \end{equation}
Then there holds for $\forall t\in(0,T]$
\begin{align}\label{4.4}
&\|u\|_{\tilde{L}^\infty_t(\dot{B}_{2,1}^{1/2})}+
\|u\|_{L^1_t(\dot{B}_{2,1}^{5/2})}+\|\nabla P\|_{L^1_t(\dot{B}_{2,1}^{1/2})}\nonumber\\
&\le C\Bigg(\|u_0\|_{\dot{B}_{2,1}^{1/2}}+\|f\|_{L^1_t(\dot{B}_{2,1}^{1/2})}+2^{2m}\|\dot{S}_ma\|_{L_t^\infty(L^2)}\|\mathrm{div }f\|_{L_t^1(H^{-1})}\nonumber\\
&~~~+\int_0^t\|u(\tau)\|_{\dot{B}_{2,1}^{1/2}}\|v(\tau)\|_{\dot{B}_{2,1}^{5/2}}d\tau+W(t)\int_0^t\|u(\tau)\|_{\dot{B}_{2,1}^{1/2}}d\tau
\Bigg),
\end{align}
with
$$W(t)=2^{8m}\|\dot{S}_ma\|^4_{L_t^\infty(L^2)}(1+\|a\|^4_{L_t^\infty(L^\infty)}+\|v\|^4_{L^\infty_t(\dot{B}_{2,1}^{1/2})})
+2^{2m}\|a\|^2_{L_t^\infty(L^\infty)},$$
provided that
$$\|a-\dot{S}_ma\|_{L^\infty_T(\dot{B}_{2,1}^{3/2})}\le c$$
for some sufficiently small positive constant $c$ and some integer $m\in\Z$.
\end{proposition}

Applying Proposition \ref{dongliang} to the second equation of (\ref{feixianxing}), we have
\begin{align}\label{4.5}
&\|\bar{u}^n\|_{\tilde{L}^\infty_t(\dot{B}_{2,1}^{1/2})}+
\|\bar{u}^n\|_{L^1_t(\dot{B}_{2,1}^{5/2})}+\|\nabla P^n\|_{L^1_t(\dot{B}_{2,1}^{1/2})}\nonumber\\
\le& C\Bigg(\|M_n\|_{L^1_t(\dot{B}_{2,1}^{1/2})}+2^{2m}\|\dot{S}_ma^n\|_{L_t^\infty(L^2)}\|\mathrm{div } M_n\|_{L_t^1(\dot{H}^{-1})}\nonumber\\
&+\int_0^t\|\bar{u}^n(\tau)\|_{\dot{B}_{2,1}^{1/2}}\|u_L^n(\tau)\|_{\dot{B}_{2,1}^{5/2}}d\tau+W^n(t)\int_0^t\|\bar{u}^n(\tau)\|_{\dot{B}_{2,1}^{1/2}}d\tau
\Bigg),
\end{align}
with
$$W^n(t)=2^{8m}\|\dot{S}_ma^n\|^4_{L_t^\infty(L^2)}(1+\|a^n\|^4_{L_t^\infty(L^\infty)}+\|u_L^n\|^4_{L^\infty_t(\dot{B}_{2,1}^{1/2})})
+2^{2m}\|a^n\|^2_{L_t^\infty(L^\infty)},$$
\begin{align*}
M_n(\tau)=&-u_L^n\cdot \nabla u_L^n-\bar{u}^n\cdot \nabla u_L^n-\bar{u}^n\cdot \nabla \bar{u}^n+a^n\Delta u_L^n
\\&+(1+a^n)(\nabla d_L^n\cdot \Delta d_L^n+\nabla d_L^n\cdot \Delta \bar{d}^n+\nabla \bar{d}^n\cdot \Delta \bar{d}^n+\nabla \bar{d}^n\cdot \Delta d_L^n),
\end{align*}
under the assumption
$$\|a^n-\dot{S}_ma^n\|_{L^\infty_T(\dot{B}_{2,1}^{3/2})}\le c.$$
Noticing that $\|fg\|_{\dot{H}^{-1}}\le C \|f\|_{\dot{B}_{2,1}^{1/2}}\|g\|_{\dot{B}_{2,1}^0}$ and
$$ \|f\cdot\nabla g\|_{L^2}\le C \|f\|_{L^3}\|\nabla g\|_{L^6}\le C \|f\|_{\dot{B}_{2,1}^{1/2}}\|g\|^{1/4}_{\dot{B}_{2,1}^{1/2}}\|g\|^{3/4}_{\dot{B}_{2,1}^{5/2}},$$
we have
\begin{align}\label{4.6}
&\|\mathrm{div}(-u_L^n\cdot \nabla u_L^n-\bar{u}^n\cdot \nabla u_L^n-\bar{u}^n\cdot \nabla \bar{u}^n+a^n\Delta u_L^n)\|_{\dot{H}^{-1}}\nonumber\\
\lesssim& \|u_L^n\|^{5/4}_{\dot{B}_{2,1}^{1/2}}\|u_L^n\|^{3/4}_{\dot{B}_{2,1}^{5/2}}+\|\bar{u}^n\|_{\dot{B}_{2,1}^{1/2}}
\|u_L^n\|^{1/4}_{\dot{B}_{2,1}^{1/2}}\|u_L^n\|^{3/4}_{\dot{B}_{2,1}^{5/2}}
+\|\bar{u}^n\|^{5/4}_{\dot{B}_{2,1}^{1/2}}\|\bar{u}^n\|^{3/4}_{\dot{B}_{2,1}^{5/2}}\nonumber\\
&+\|a^n\|_{\dot{B}_{2,1}^{3/2}}\|u_L^n\|^{1/4}_{\dot{B}_{2,1}^{1/2}}\|u_L^n\|^{3/4}_{\dot{B}_{2,1}^{5/2}}
\end{align}
and
\begin{align}\label{4.7}
&\|\mathrm{div}((1+a^n)(\nabla d_L^n\cdot \Delta d_L^n+\nabla d_L^n\cdot \Delta \bar{d}^n+\nabla \bar{d}^n\cdot \Delta \bar{d}^n+\nabla \bar{d}^n\cdot \Delta d_L^n))\|_{\dot{H}^{-1}}\nonumber\\
\lesssim & (1+\|a^n\|_{\dot{B}_{2,1}^{3/2}})(\|d_L^n\|^{5/4}_{\dot{B}_{2,1}^{3/2}}\|d_L^n\|^{3/4}_{\dot{B}_{2,1}^{7/2}}
+\|d_L^n\|_{\dot{B}_{2,1}^{3/2}}\|\bar{d}^n\|^{1/4}_{\dot{B}_{2,1}^{3/2}}\|\bar{d}^n\|^{3/4}_{\dot{B}_{2,1}^{7/2}}\nonumber\\
&+\|\bar{d}^n\|_{\dot{B}_{2,1}^{3/2}}\|d_L^n\|^{1/4}_{\dot{B}_{2,1}^{3/2}}\|d_L^n\|^{3/4}_{\dot{B}_{2,1}^{7/2}}
+\|\bar{d}^n\|^{5/4}_{\dot{B}_{2,1}^{3/2}}\|\bar{d}^n\|^{3/4}_{\dot{B}_{2,1}^{7/2}}).
\end{align}
Combining with (\ref{4.6}) and (\ref{4.7}) and applying the Young's inequality yields that
\begin{align}\label{4.8}
&2^{2m}\|\dot{S}_ma^n\|_{L_t^\infty(L^2)}\|\mathrm{div}M_n\|_{L^1_t(\dot{H}^{-1})}\nonumber\\
\le & \varepsilon(\|\bar{u}^n\|_{L_t^1({\dot{B}_{2,1}^{5/2}})}+\|\bar{d}^n\|_{L_t^1({\dot{B}_{2,1}^{7/2}})})
+\|u_L^n\|_{L_t^1({\dot{B}_{2,1}^{5/2}})}+\|d_L^n\|_{L_t^1({\dot{B}_{2,1}^{7/2}})}\nonumber\\
&+Ct2^{8m}\|\dot{S}_ma^n\|^4_{L_t^\infty(L^2)}(1+\|a^n\|^4_{\tilde{L}^\infty_t(\dot{B}_{2,1}^{3/2})})  (Z^n(t))^5\nonumber\\
&+Ct
2^{8m}\|\dot{S}_ma^n\|^4_{L_t^\infty(L^2)}(1+\|a^n\|^4_{\tilde{L}^\infty_t(\dot{B}_{2,1}^{3/2})})
(1+\|u_0\|_{\dot{B}_{2,1}^{1/2}}^5+\|d_0\|_{\dot{B}_{2,1}^{3/2}}^5)
\nonumber\\
&+Ct2^{8m}\|\dot{S}_ma^n\|^4_{L_t^\infty(L^2)}(1+\|a^n\|^4_{\tilde{L}^\infty_t(\dot{B}_{2,1}^{3/2})})
(\|u_0\|_{\dot{B}_{2,1}^{1/2}}+\|d_0\|_{\dot{B}_{2,1}^{3/2}})(Z^n(t))^4
\nonumber\\
&+Ct2^{8m}\|\dot{S}_ma^n\|^4_{L_t^\infty(L^2)}(1+\|a^n\|^4_{\tilde{L}^\infty_t(\dot{B}_{2,1}^{3/2})})
\|d_0\|_{\dot{B}_{2,1}^{3/2}}^4(Z^n(t)).
\end{align}
Applying product laws in Besov spaces gives

\begin{align}\label{4.9}
\|M_n\|_{L_t^1({\dot{B}_{2,1}^{1/2}})}\lesssim& \int_0^t\|u_L^n\|_{\dot{B}_{2,1}^{1/2}}\|u_L^n\|_{\dot{B}_{2,1}^{5/2}}+\|\bar{u}^n\|_{\dot{B}_{2,1}^{1/2}}\|\bar{u}^n\|_{\dot{B}_{2,1}^{5/2}}
+\|\bar{u}^n\|_{\dot{B}_{2,1}^{1/2}}\|u_L^n\|_{\dot{B}_{2,1}^{5/2}}d\tau\nonumber\\
&+(1+\|a^n\|_{\tilde{L}^\infty_t(\dot{B}_{2,1}^{3/2})}) \int_0^t\Big(\|u_L^n\|_{\dot{B}_{2,1}^{5/2}}+\|\nabla d_L^n\|_{\dot{B}_{2,1}^{3/2}}\|\Delta d_L^n\|_{\dot{B}_{2,1}^{1/2}}+\|\nabla d_L^n\|_{\dot{B}_{2,1}^{3/2}}\|\Delta \bar{d}^n\|_{\dot{B}_{2,1}^{1/2}}\nonumber\\
&\hspace{3.5cm} +\|\nabla \bar{d}^n\|_{\dot{B}_{2,1}^{3/2}}\|\Delta \bar{d}^n\|_{\dot{B}_{2,1}^{1/2}}+\|\nabla \bar{d}^n\|_{\dot{B}_{2,1}^{3/2}}\|\Delta d_L^n\|_{\dot{B}_{2,1}^{1/2}}\Big)d\tau\nonumber\\
\lesssim& \int_0^t\|u_L^n\|_{\dot{B}_{2,1}^{1/2}}\|u_L^n\|_{\dot{B}_{2,1}^{5/2}}+\|\bar{u}^n\|_{\dot{B}_{2,1}^{1/2}}\|\bar{u}^n\|_{\dot{B}_{2,1}^{5/2}}
+\|\bar{u}^n\|_{\dot{B}_{2,1}^{1/2}}\|u_L^n\|_{\dot{B}_{2,1}^{5/2}}d\tau\nonumber\\
&+(1+\|a^n\|_{\tilde{L}^\infty_t(\dot{B}_{2,1}^{3/2})}) \int_0^t\|u_L^n\|_{\dot{B}_{2,1}^{5/2}}+\| d_L^n\|_{\dot{B}_{2,1}^{5/2}}^2+\|d_L^n\|_{\dot{B}_{2,1}^{5/2}}\| \bar{d}^n\|_{\dot{B}_{2,1}^{5/2}}
+\|\bar{d}^n\|_{\dot{B}_{2,1}^{5/2}}^2d\tau\nonumber\\
\lesssim& \int_0^t\|u_L^n\|_{\dot{B}_{2,1}^{1/2}}\|u_L^n\|_{\dot{B}_{2,1}^{5/2}}+\|\bar{u}^n\|_{\dot{B}_{2,1}^{1/2}}\|\bar{u}^n\|_{\dot{B}_{2,1}^{5/2}}
+\|\bar{u}^n\|_{\dot{B}_{2,1}^{1/2}}\|u_L^n\|_{\dot{B}_{2,1}^{5/2}}d\tau\nonumber\\
&+(1+\|a^n\|_{\tilde{L}^\infty_t(\dot{B}_{2,1}^{3/2})}) \int_0^t\Big(\|u_L^n\|_{\dot{B}_{2,1}^{5/2}}+\| d_L^n\|_{\dot{B}_{2,1}^{3/2}}\| d_L^n\|_{\dot{B}_{2,1}^{7/2}}
+\| \bar{d}^n\|_{\dot{B}_{2,1}^{3/2}}\| \bar{d}^n\|_{\dot{B}_{2,1}^{7/2}}
\Big)d\tau\nonumber\\
\lesssim& (1+\|a^n\|_{\tilde{L}^\infty_t(\dot{B}_{2,1}^{3/2})})( \|u_0\|_{\dot{B}_{2,1}^{1/2}}+\|d_0\|_{\dot{B}_{2,1}^{3/2}})
(\|u_L^n\|_{L^1_t(\dot{B}_{2,1}^{5/2})}+\|d_L^n\|_{L^1_t(\dot{B}_{2,1}^{7/2})})\nonumber\\
&+(Z^n(t))^2+\int_0^t\|u_L^n\|_{\dot{B}_{2,1}^{5/2}}Z^n(\tau)d\tau.
\end{align}
Now, we define $m\in\Z$ by
$$m=\inf\left\{p\in\Z|\sum_{q\ge p}2^{3q/2} \|\Delta_qa_0\|_{L^2}\le c_0\underline{b}\right\}$$
for some sufficiently small positive constant $c_0$.

Noticing that $\mathrm{div}( \bar{u}^n +u_L^n) = 0,$ we get by applying Proposition \ref{shuyun} to the first
equation of (\ref{feixianxing}) that
\begin{align}\label{4.10}
\|a^n\|_{\tilde{L}^\infty_t(B_{2,1}^{3/2})}\le\|a_0\|_{B_{2,1}^{3/2}}\exp\left\{ C \(\|u_L^n\|_{L_t^1({\dot{B}_{2,1}^{5/2}})}+\|\bar{u}^n\|_{L_t^1({\dot{B}_{2,1}^{5/2}})}\)\right\}
\end{align}
and
\begin{align}\label{4.11}
\|a^n\|_{L^\infty_t(L^2)}\le\|a_0\|_{L^2}, \hspace{0.5cm}\|a^n\|_{L^\infty_t(L^\infty)}\le\|a_0\|_{L^\infty},
\end{align}
from which, we can get
\begin{align}\label{4.12}
W^n(t)=&2^{8m}\|\dot{S}_ma^n\|^4_{L_t^\infty(L^2)}(1+\|a^n\|^4_{L_t^\infty(L^\infty)}+\|u_L^n\|^4_{L^\infty_t(\dot{B}_{2,1}^{1/2})})
+2^{2m}\|a^n\|^2_{L_t^\infty(L^\infty)}\nonumber\\
&\le2^{8m}\|a_0\|^4_{L^2}\(1+\|a_0\|_{L^\infty}^4+\|u_0\|^4_{\dot{B}_{2,1}^{1/2}}\)+2^{2m}\|a_0\|^2_{L^\infty} \triangleq N_m.
\end{align}
Inserting (\ref{4.8}) (\ref{4.9}) into (\ref{4.5}) yields that
\begin{align}\label{4.13}
&\|\bar{u}^n\|_{\tilde{L}^\infty_t(\dot{B}_{2,1}^{1/2})}+
\|\bar{u}^n\|_{L^1_t(\dot{B}_{2,1}^{5/2})}+\|\nabla P^n\|_{L^1_t(\dot{B}_{2,1}^{1/2})}\nonumber\\
\le&\varepsilon(\|\bar{u}^n\|_{L_t^1({\dot{B}_{2,1}^{5/2}})}+\|\bar{d}^n\|_{L_t^1({\dot{B}_{2,1}^{7/2}})})
 +C\int_0^t(N_m+\|u_L^n\|_{\dot{B}_{2,1}^{5/2}}+\|d_L^n\|_{\dot{B}_{2,1}^{7/2}})Z^n(\tau)d\tau\nonumber\\
&+C(Z^n(t))^2+C(1+\|a^n\|_{\tilde{L}^\infty_t(\dot{B}_{2,1}^{3/2})})( 1+ \|u_0\|_{\dot{B}_{2,1}^{1/2}}+\|d_0\|_{\dot{B}_{2,1}^{3/2}})
(\|u_L^n\|_{L^1_t(\dot{B}_{2,1}^{5/2})}+\|d_L^n\|_{L^1_t(\dot{B}_{2,1}^{7/2})})\nonumber\\
&+Ct2^{8m}\|\dot{S}_ma^n\|^4_{L_t^\infty(L^2)}(1+\|a^n\|^4_{\tilde{L}^\infty_t(\dot{B}_{2,1}^{3/2})})  (Z^n(t))^5\nonumber\\
&+Ct
2^{8m}\|\dot{S}_ma^n\|^4_{L_t^\infty(L^2)}(1+\|a^n\|^4_{\tilde{L}^\infty_t(\dot{B}_{2,1}^{3/2})})
(1+\|u_0\|_{\dot{B}_{2,1}^{1/2}}^5+\|d_0\|_{\dot{B}_{2,1}^{3/2}}^5)
\nonumber\\
&+Ct2^{8m}\|\dot{S}_ma^n\|^4_{L_t^\infty(L^2)}(1+\|a^n\|^4_{\tilde{L}^\infty_t(\dot{B}_{2,1}^{3/2})})
(\|u_0\|_{\dot{B}_{2,1}^{1/2}}+\|d_0\|_{\dot{B}_{2,1}^{3/2}})(Z^n(t))^4
\nonumber\\
&+Ct2^{8m}\|\dot{S}_ma^n\|^4_{L_t^\infty(L^2)}(1+\|a^n\|^4_{\tilde{L}^\infty_t(\dot{B}_{2,1}^{3/2})})
\|d_0\|_{\dot{B}_{2,1}^{3/2}}^4(Z^n(t)).
\end{align}
Summing the above result  up  to (\ref{d11}), together with (\ref{4.1}), (\ref{4.10}), (\ref{4.11}) and taking $\varepsilon$ small enough, we can finally get
\begin{align}\label{4.15}
Z^n(t)
=&\|\bar{u}^n\|_{\tilde{L}^\infty_t(\dot{B}_{2,1}^{1/2})}+\|\bar{d}^n\|_{\tilde{L}^\infty_t(\dot{B}_{2,1}^{3/2})}
+\|\bar{u}^n\|_{L_t^1({\dot{B}_{2,1}^{5/2}})}+\|\bar{d}^n\|_{L_t^1({\dot{B}_{2,1}^{7/2}})}
+\|\nabla P^n\|_{L_t^1({\dot{B}_{2,1}^{1/2}})}\nonumber\\
\lesssim&
 \int_0^t(N_m+\|u_L^n\|_{\dot{B}_{2,1}^{5/2}}+(1+\|d_0\|_{\dot{B}_{2,1}^{3/2}})\|d_L^n\|_{\dot{B}_{2,1}^{7/2}})Z^n(\tau)d\tau
+C(1+\|d_0\|_{\dot{B}_{2,1}^{3/2}}+(Z^n(t))^2)(Z^n(t))^2\nonumber\\
&+(1+\|a_0\|_{\dot{B}_{2,1}^{3/2}}\exp(C\|u_0\|_{\dot{B}_{2,1}^{1/2}})e^{(Z^n(t))})( 1+ \|u_0\|_{\dot{B}_{2,1}^{1/2}}+\|d_0\|_{\dot{B}_{2,1}^{3/2}})
(\|u_L^n\|_{L^1_t(\dot{B}_{2,1}^{5/2})}+\|d_L^n\|_{L^1_t(\dot{B}_{2,1}^{7/2})})\nonumber\\
&+t2^{8m}\|a_0\|_{L^2}^4 (1+\|a_0\|^4_{\dot{B}_{2,1}^{3/2}}\exp(C\|u_0\|_{\dot{B}_{2,1}^{1/2}})e^{(Z^n(t))})
(\|u_0\|_{\dot{B}_{2,1}^{1/2}}+\|d_0\|_{\dot{B}_{2,1}^{3/2}}+Z^n(t))(Z^n(t))^4
\nonumber\\
&+t
2^{8m} \|a_0\|_{L^2}^4 (1+\|a_0\|^4_{\dot{B}_{2,1}^{3/2}}\exp(C\|u_0\|_{\dot{B}_{2,1}^{1/2}})e^{(Z^n(t))})
(1+\|u_0\|_{\dot{B}_{2,1}^{1/2}}^5+\|d_0\|_{\dot{B}_{2,1}^{3/2}}^5)
\nonumber\\
&+t2^{8m}\|a_0\|_{L^2}^4 (1+\|a_0\|^4_{\dot{B}_{2,1}^{3/2}}\exp(C\|u_0\|_{\dot{B}_{2,1}^{1/2}})e^{(Z^n(t))})
\|d_0\|_{\dot{B}_{2,1}^{3/2}}^4(Z^n(t)),
\end{align}
under the assumption
$$\|a^n-\dot{S}_ma^n\|_{L^\infty_t(\dot{B}_{2,1}^{3/2})}\le 2c_0\underline{b}.$$

Gronwall inequality help us to get
\begin{align}\label{4.19}
Z^n(t)\le& C_3\exp(C_3(tN_m+\|u_0\|_{\dot{B}_{2,1}^{1/2}}+(1+\|d_0\|_{\dot{B}_{2,1}^{3/2}})\|d_0\|_{\dot{B}_{2,1}^{3/2}}   ))\nonumber\\
&\times\Bigg\{C_3(1+\|d_0\|_{\dot{B}_{2,1}^{3/2}}+Z^n(t))(Z^n(t))^2\nonumber\\
&+C_3(1+\|a_0\|_{\dot{B}_{2,1}^{3/2}}\exp(C_3\|u_0\|_{\dot{B}_{2,1}^{1/2}})e^{(Z^n(t))})( 1+ \|u_0\|_{\dot{B}_{2,1}^{1/2}}+\|d_0\|_{\dot{B}_{2,1}^{3/2}})\nonumber\\
&\times
\(\sum_{q\in\Z}(1-e^{-ct2^{2q}})
(2^{q/2}\|\dot{\Delta}_qu_0\|_{L^2}+2^{{3q}/2}\|\dot{\Delta}_qd_0\|_{L^2})\)\nonumber\\
&+C_3t2^{8m}\|a_0\|_{L^2}^4 (1+\|a_0\|^4_{\dot{B}_{2,1}^{3/2}}\exp(C_3\|u_0\|_{\dot{B}_{2,1}^{1/2}})e^{(Z^n(t))})
(\|u_0\|_{\dot{B}_{2,1}^{1/2}}+\|d_0\|_{\dot{B}_{2,1}^{3/2}}+Z^n(t))(Z^n(t))^4
\nonumber\\
&+C_3t
2^{8m} \|a_0\|_{L^2}^4 (1+\|a_0\|^4_{\dot{B}_{2,1}^{3/2}}\exp(C_3\|u_0\|_{\dot{B}_{2,1}^{1/2}})e^{(Z^n(t))})
(1+\|u_0\|_{\dot{B}_{2,1}^{1/2}}^5+\|d_0\|_{\dot{B}_{2,1}^{3/2}}^5)
\nonumber\\
&+C_3t2^{8m}\|a_0\|_{L^2}^4 (1+\|a_0\|^4_{\dot{B}_{2,1}^{3/2}}\exp(C_3\|u_0\|_{\dot{B}_{2,1}^{1/2}})e^{(Z^n(t))})
\|d_0\|_{\dot{B}_{2,1}^{3/2}}^4(Z^n(t))\Bigg\},
\end{align}
where we assume the constant $C_3>C_2>1.$

Applying (\ref{3.3}) to the first equation of (\ref{feixianxing}) together with (\ref{4.0}),(\ref{4.1}), (\ref{4.10}), (\ref{4.11}) and $1+xe^x\ge e^x$ for $\forall x\ge 0$, we deduce that
\begin{align}\label{4.16}
&\|a^n-\dot{S}_ma^n\|_{\tilde{L}^\infty_t(\dot{B}_{2,1}^{3/2})}\le \|a^n-S_ma^n\|_{\tilde{L}^\infty_t(B_{2,1}^{3/2})}\nonumber\\
\le&\sum_{q\ge m} 2^{3q/2}\|\Delta_qa_0\|_{L^2}+\|a_0\|_{B_{2,1}^{3/2}}\(\exp{C(Z^n(t)+\|u_L^n\|_{L_t^1({\dot{B}_{2,1}^{5/2}})})}-1\)\nonumber\\
\le &c_0\underline{b}+\|a_0\|_{B_{2,1}^{3/2}} (Z^n(t)+\|u_L^n\|_{L_t^1({\dot{B}_{2,1}^{5/2}})})\exp\(C(Z^n(t)+\|u_L^n\|_{L_t^1({\dot{B}_{2,1}^{5/2}})})\)\nonumber\\
\le& c_0\underline{b}+C_2\|a_0\|_{B_{2,1}^{3/2}}\exp(C_2\|u_0\|_{\dot{B}_{2,1}^{1/2}})e^{(C_2Z^n(t))}\(\sum_{q\in\Z}2^{q/2}(1-e^{-ct2^{2q}})
\|\dot{\Delta}_qu_0\|_{L^2}+Z^n(t)\).
\end{align}
Using the fact that $1+x\le e^x,$ $x\ge 0$ and  taking $0<T_1(\varepsilon_0)\le(C_3N_m)^{-1}$ so small that
\begin{align}\label{4.20}
eC_3(1+\|a_0\|_{B_{2,1}^{3/2}})&\exp(C_3)\nonumber\\
&\times\(\sum_{q\in\Z}(1-e^{-c_1T_12^{2q}})
(2^{q/2}\|\dot{\Delta}_qu_0\|_{L^2}+2^{{3q}/2}\|\dot{\Delta}_qd_0\|_{L^2})\)\le\varepsilon_0,
\end{align}
for some sufficiently small positive constant $\varepsilon_0$.

Then it follows from (\ref{4.19})
that
\begin{align}\label{4.21}
Z^n(t)\le&\varepsilon_0e^{C_3Z^n(t)}+ eC_3\exp(C_3(\|u_0\|_{\dot{B}_{2,1}^{1/2}}+(1+\|d_0\|_{\dot{B}_{2,1}^{3/2}})\|d_0\|_{\dot{B}_{2,1}^{3/2}}))\nonumber\\
&\times\Bigg\{C_3(1+\|d_0\|_{\dot{B}_{2,1}^{3/2}}+Z^n(t))(Z^n(t))^2\nonumber\\
&+C_3t2^{8m}\|a_0\|_{L^2}^4 (1+\|a_0\|^4_{\dot{B}_{2,1}^{3/2}}\exp(C_3\|u_0\|_{\dot{B}_{2,1}^{1/2}})e^{(Z^n(t))})
(\|u_0\|_{\dot{B}_{2,1}^{1/2}}+\|d_0\|_{\dot{B}_{2,1}^{3/2}}+Z^n(t))(Z^n(t))^4
\nonumber\\
&+C_3t
2^{8m} \|a_0\|_{L^2}^4 (1+\|a_0\|^4_{\dot{B}_{2,1}^{3/2}}\exp(C_3\|u_0\|_{\dot{B}_{2,1}^{1/2}})e^{(Z^n(t))})
(1+\|u_0\|_{\dot{B}_{2,1}^{1/2}}^5+\|d_0\|_{\dot{B}_{2,1}^{3/2}}^5)
\nonumber\\
&+C_3t2^{8m}\|a_0\|_{L^2}^4 (1+\|a_0\|^4_{\dot{B}_{2,1}^{3/2}}\exp(C_3\|u_0\|_{\dot{B}_{2,1}^{1/2}})e^{(Z^n(t))})
\|d_0\|_{\dot{B}_{2,1}^{3/2}}^4(Z^n(t))\Bigg\},
\end{align}
provided that
\begin{align}\label{4.210}
C_2\|a_0\|_{B_{2,1}^{3/2}}\exp(C_2\|u_0\|_{\dot{B}_{2,1}^{1/2}})e^{(C_2Z^n(t))}\(\sum_{q\in\Z}2^{q/2}(1-e^{-cT_12^{2q}})
\|\dot{\Delta}_qu_0\|_{L^2}+Z^n(t)\)\le c_0\underline{b}.
\end{align}
Let $T^n_*(\varepsilon_0)=\sup\{t\in[0,T_1(\varepsilon_0)]|Z^n(t)\le 2e\varepsilon_0\}$. Without loss of generality, we
may assume that $\varepsilon_0$ is so small such that
$2eC_3\varepsilon_0\le1$ and $2e^2C_2\varepsilon_0\|a_0\|_{B_{2,1}^{3/2}}\exp(C_2\|u_0\|_{\dot{B}_{2,1}^{1/2}})\le c_0\underline{b}.$
Then for $t\le T^n_*(\varepsilon_0)$, we get from (\ref{4.20}) that (\ref{4.210}) holds and then
\begin{align}\label{4.22}
Z^n(t)\le&e\varepsilon_0+ eC_3\exp(C_3(\|u_0\|_{\dot{B}_{2,1}^{1/2}}+(1+\|d_0\|_{\dot{B}_{2,1}^{3/2}})\|d_0\|_{\dot{B}_{2,1}^{3/2}}))
(C_3(1+\|d_0\|_{\dot{B}_{2,1}^{3/2}}+2e\varepsilon_0)2e\varepsilon_0Z^n(t))\nonumber\\
&+eC_3\exp(C_3(\|u_0\|_{\dot{B}_{2,1}^{1/2}}+(1+\|d_0\|_{\dot{B}_{2,1}^{3/2}})\|d_0\|_{\dot{B}_{2,1}^{3/2}}))(t\mathcal{A}_1+t\mathcal{A}_2+t\mathcal{A}_3),
\end{align}
with
\begin{align}\label{}
\mathcal{A}_1=C_32^{8m}\|a_0\|_{L^2}^4 (1+e\|a_0\|^4_{\dot{B}_{2,1}^{3/2}}\exp(C_3\|u_0\|_{\dot{B}_{2,1}^{1/2}}))
(\|u_0\|_{\dot{B}_{2,1}^{1/2}}+\|d_0\|_{\dot{B}_{2,1}^{3/2}}+2e\varepsilon_0)(2e\varepsilon_0)^4,
\end{align}
\begin{align}\label{}
\mathcal{A}_2=C_3
2^{8m} \|a_0\|_{L^2}^4  (1+e\|a_0\|^4_{\dot{B}_{2,1}^{3/2}}\exp(C_3\|u_0\|_{\dot{B}_{2,1}^{1/2}}))
(1+\|u_0\|_{\dot{B}_{2,1}^{1/2}}^5+\|d_0\|_{\dot{B}_{2,1}^{3/2}}^5)
\end{align}
and
\begin{align}\label{}
\mathcal{A}_3=C_32^{8m}\|a_0\|_{L^2}^4 (1+e\|a_0\|^4_{\dot{B}_{2,1}^{3/2}}\exp(C_3\|u_0\|_{\dot{B}_{2,1}^{1/2}}))
\|d_0\|_{\dot{B}_{2,1}^{3/2}}^4(2e\varepsilon_0).
\end{align}

Taking $\varepsilon_0$ so small that
\begin{align}\label{4.23}
2e^2C_3^2\varepsilon_0(1+\|d_0\|_{\dot{B}_{2,1}^{3/2}}+2e\varepsilon_0)\exp(C_3(\|u_0\|_{\dot{B}_{2,1}^{1/2}}+(1+\|d_0\|_{\dot{B}_{2,1}^{3/2}})\|d_0\|_{\dot{B}_{2,1}^{3/2}}))
\le{1}/{15}
\end{align}
and
\begin{align}\label{4.24}
T_2(\varepsilon_0)=\min\(T_1(\varepsilon_0),t_1,t_2,t_3\)
\end{align}
where
\begin{align}\label{4.25}
t_1\le{e\varepsilon_0}/{\left(eC_3^2\mathcal{A}_1\exp(C_3(\|u_0\|_{\dot{B}_{2,1}^{1/2}}+(1+\|d_0\|_{\dot{B}_{2,1}^{3/2}})\|d_0\|_{\dot{B}_{2,1}^{3/2}}))\right)},
\end{align}
\begin{align}\label{4.26}
t_2\le{e\varepsilon_0}/{\left(eC_3^2\mathcal{A}_2\exp(C_3(\|u_0\|_{\dot{B}_{2,1}^{1/2}}+(1+\|d_0\|_{\dot{B}_{2,1}^{3/2}})\|d_0\|_{\dot{B}_{2,1}^{3/2}}))\right)}
\end{align}
and
\begin{align}\label{4.27}
t_3\le{e\varepsilon_0}/{\left(eC_3^2\mathcal{A}_3\exp(C_3(\|u_0\|_{\dot{B}_{2,1}^{1/2}}+(1+\|d_0\|_{\dot{B}_{2,1}^{3/2}})\|d_0\|_{\dot{B}_{2,1}^{3/2}}))\right)}.
\end{align}
Therefore, we have
\begin{align}\label{4.28}
Z^n(t)\le{4e\varepsilon_0}/{3},\hspace{0.5cm} for \ \forall\le T_2(\varepsilon_0).
\end{align}
This, together with (\ref{4.0}) and (\ref{4.1}), ensures
\begin{align}\label{4.29}
\{a^n,u^n,d^n,\nabla P^n\}_{n\in\N} \,\,\,\,{\rm is}\,\,{\rm uniformly}\,\,{\rm bounded}\,\,{\rm in}\,\, E_{T_2(\varepsilon_0)}.
\end{align}

STEP 3: (Convergence).
 The convergence of $(u_L^n,  d_L^n)$ to $(u_L, d_L)$ readily stems from the definition of Besov spaces.
As for the convergence of $(a^n,\bar{u}^n,\bar{d}^n)$, it relies upon Ascoli's theorem compactness properties of the consequence, which are obtained by considering the time derivative of the solution, we omit the details here .Moreover, there holds
\begin{align}\label{xiaotiaojian}
\|a^n-S_ma^n\|_{L^\infty_t(B_{2,1}^{3/2})}\le 2c_0\underline{b}
\end{align}
with the constant $c_0$ being sufficiently small.

\subsection{Uniqueness of Theorem 1.1 }
Let $(a^i , u^i , d^i,\nabla P^i )$ (with $i = 1, 2$) be two solutions of the system  (\ref{1.2}), which
satisfy  (\ref{xiaotiaojian}) and
\begin{align}\label{5.0}
(a^i, u^i,d^i,\nabla P^i ) \in & C_b([0, T ]; B_{2,1}^{3/2}\times C_b([0, T ]; B_{2,1}^{1/2}\cap L^1([0,T];B_{2,1}^{5/2})\nonumber\\
&\times C_b([0, T ]; B_{2,1}^{3/2} )
\cap L^1([0,T];B_{2,1}^{7/2})\times L^1([0,T];B_{2,1}^{1/2}).
\end{align}
We define
$$ (\delta a,\delta u,\delta d,\nabla\delta P)=(a^2-a^1,u^2-u^1,d^2-d^1,\nabla P^2-\nabla P^1,$$
so that $(\delta a,\delta u,\delta d,\nabla\delta P)$ solves
\begin{eqnarray}\label{5.1}
\left\{\begin{aligned}
&\partial_t\delta a+u^2\cdot\nabla \delta a=-\delta u\cdot\nabla a^1,\\
&\partial_t\delta u+u^2 \cdot \nabla \delta u+(1+a^2)(\nabla\delta P-\Delta \delta u)=F,\\
&\partial_t\delta d-\Delta \delta d+u^2 \cdot \nabla \delta d=\delta u\cdot \nabla d^1+\nabla d^2\cdot\nabla\delta d\cdot d^2+\nabla \delta d\cdot\nabla d^1\cdot d^2+\nabla d^1\cdot\nabla d^1\cdot \delta d,\\
&\mathrm{\mathrm{div}}\delta u=0,\\
&(\delta a,\delta u,\delta d)|_{t=0}=(0,0,0),
\end{aligned}\right.
\end{eqnarray}
with

\begin{align*}
F=& -\delta u\cdot \nabla u^1+\delta a(\Delta u^1-\nabla P^1)
+(1+a^1)(\nabla d^1\cdot\Delta \delta d)+(1+a^2)(\nabla \delta d \cdot\Delta d^2)+\delta a \nabla d^1\cdot\Delta d^2.
\end{align*}

Applying Proposition \ref{3.0} to the first equation in (\ref{5.1}) yields
\begin{align}\label{5.2}
\|\delta a\|_{\tilde{L}^\infty_t(B_{2,1}^{1/2})}\le\exp\{C\|\nabla u^2\|_{L_t^1({B_{2,1}^{3/2}})}\}\|a^1\|_{\tilde{L}^\infty_t(B_{2,1}^{3/2})}\|\delta u \|_{L^1_t(B_{2,1}^{3/2})}.
\end{align}

For the momentum equation of (\ref{feixianxing}), we can follow the proof of \cite{abidi2012} (up to a slight modification) to get
\begin{align}\label{5.12}
&\|\delta u\|_{\tilde{L}^\infty_t(B_{2,1}^{-1/2})}+
\|\delta u\|_{L^1_t(B_{2,1}^{3/2})}\nonumber\\
\le&   C\int_0^t\|\delta u\|_{B_{2,1}^{-1/2}}(1+\|\nabla u^1\|_{B_{2,1}^{3/2}}^2+\|\nabla u^2\|_{B_{2,1}^{3/2}}^2+\| u^2\|_{B_{2,1}^{3/2}}^2)d\tau+C\int_0^t\|\delta a\|_{B_{2,1}^{1/2}}\| d^1\|_{B_{2,1}^{3/2}}\| d^2\|_{B_{2,1}^{7/2}}d\tau\nonumber\\
&+   C\int_0^t\|\delta d\|_{B_{2,1}^{1/2}}(1+\| a^1\|^2_{B_{2,1}^{3/2}}+\| a^2\|_{B_{2,1}^{3/2}})(\| d^1\|_{B_{2,1}^{3/2}}+\| d^2\|_{B_{2,1}^{3/2}})(1+\| d^1\|_{B_{2,1}^{7/2}}+\| d^2\|_{B_{2,1}^{7/2}})d\tau\nonumber\\
&+C\int_0^t\|\delta u\|_{L_\tau^1({B_{2,1}^{3/2}})}(\|\Delta u^1\|_{B_{2,1}^{1/2}}+\|\nabla P^1\|_{B_{2,1}^{1/2}}+\|d^1\|_{B_{2,1}^{3/2}}\|d^2\|_{B_{2,1}^{7/2}})d\tau+\varepsilon\|\delta d\|_{L^1_t(B_{2,1}^{7/2})},
\end{align}
where we have used the following key estimates which are different to Navier-Stokes equations:
\begin{align}\label{}
\|(1+a^1)(\nabla d^1\cdot\Delta \delta d)\|_{L^1_t(B_{2,1}^{-1/2})}\le& C\int_0^t (1+\|a^1\|_{B_{2,1}^{3/2}})\|\nabla d^1\|_{B_{2,1}^{3/2}}\|\Delta \delta d\|_{B_{2,1}^{-1/2}} d\tau\nonumber\\
\le& C\int_0^t (1+\|a^1\|_{B_{2,1}^{3/2}})\| d^1\|_{B_{2,1}^{5/2}}\|\delta d\|_{B_{2,1}^{3/2}} d\tau\nonumber\\
\le& C\int_0^t (1+\|a^1\|_{B_{2,1}^{3/2}})\| d^1\|_{B_{2,1}^{3/2}}^{1/2}\| d^1\|_{B_{2,1}^{7/2}}^{1/2}\|\delta d\|_{B_{2,1}^{1/2}}^{1/2} \|\delta d\|_{B_{2,1}^{5/2}}^{1/2} d\tau\nonumber\\
\le& \varepsilon  \|\delta d\|_{L^1_t(B_{2,1}^{5/2})}+ C\int_0^t (1+\|a^1\|^2_{B_{2,1}^{3/2}})\| d^1\|_{B_{2,1}^{3/2}}\| d^1\|_{B_{2,1}^{7/2}}\|\delta d\|_{B_{2,1}^{1/2}}d\tau,
\end{align}

\begin{align}\label{}
\|(1+a^1)(\nabla d^1\cdot\Delta \delta d)\|_{L^1_t(H^{-1})}\le& C\int_0^t (1+\|a^1\|_{B_{2,1}^{3/2}})\|\nabla d^1\|_{B_{2,1}^{3/2}}\|\Delta \delta d\|_{H^{-1}} d\tau\nonumber\\
\le& C\int_0^t (1+\|a^1\|_{B_{2,1}^{3/2}})\|\nabla d^1\|_{B_{2,1}^{3/2}}\|\delta d\|_{H^{1}} d\tau\nonumber\\
\le& C\int_0^t (1+\|a^1\|_{B_{2,1}^{3/2}})\|\nabla d^1\|_{B_{2,1}^{3/2}}\|\delta d\|_{B_{2,1}^{1/2}}^{3/4}
\|\delta d\|_{B_{2,1}^{5/2}}^{1/4} d\tau\nonumber\\
\le& \varepsilon \|\delta d\|_{L^1_t(B_{2,1}^{5/2})} + C\int_0^t (1+\|a^1\|^2_{B_{2,1}^{3/2}})(1+\|d^1\|_{B_{2,1}^{5/2}}^2)\|\delta d\|_{B_{2,1}^{1/2}}d\tau\nonumber\\
\le& \varepsilon \|\delta d\|_{L^1_t(B_{2,1}^{5/2})} + C\int_0^t (1+\|a^1\|^2_{B_{2,1}^{3/2}})(1+\|d^1\|_{B_{2,1}^{3/2}}\|d^1\|_{B_{2,1}^{7/2}})\|\delta d\|_{B_{2,1}^{1/2}}d\tau
\end{align}
and
\begin{align}\label{}
\|(1+a^2)(\nabla \delta d\cdot \Delta  d^2)\|_{L^1_t(H^{-1})}\le& C\int_0^t (1+\|a^2\|_{B_{2,1}^{3/2}})\| d^2\|_{B_{2,1}^{3/2}}\|\delta d\|_{H^{1}} d\tau\nonumber\\
\le& \varepsilon \|\delta d\|_{L^1_t(B_{2,1}^{5/2})} + C\int_0^t (1+\|a^2\|^2_{B_{2,1}^{3/2}})(1+\|d^2\|_{B_{2,1}^{5/2}}^2)\|\delta d\|_{B_{2,1}^{1/2}}d\tau\nonumber\\
\le& \varepsilon \|\delta d\|_{L^1_t(B_{2,1}^{5/2})} + C\int_0^t (1+\|a^2\|^2_{B_{2,1}^{3/2}})(1+\|d^2\|_{B_{2,1}^{3/2}}\|d^2\|_{B_{2,1}^{7/2}})\|\delta , d\|_{B_{2,1}^{1/2}}d\tau,
\end{align}
\begin{align}\label{}
\|\delta a \nabla d^1\cdot\Delta d^2\|_{L^1_t(H^{-1})}\le& C\int_0^t \|\delta a\|_{B_{2,1}^{1/2}}
\|\nabla d^1\|_{L^2}\|\Delta d^2\|_{L^\infty}d\tau\nonumber\\
\le& C\int_0^t \|\delta a\|_{B_{2,1}^{1/2}}
\|d^1\|_{B_{2,1}^{3/2}}\|d^2\|_{B_{2,1}^{3/2}}d\tau.
\end{align}

Similar to the estimate (\ref{4.3}), we can get from the  third equation of (\ref{5.1}) that
\begin{align}\label{5.3}
&\|\delta d\|_{\tilde{L}^\infty_t(B_{2,1}^{1/2})}+
\|\delta d\|_{L^1_t(B_{2,1}^{5/2})}\nonumber\\
\lesssim& \sum_{q\in\Z}2^{{q}/2} \|[u^2 \cdot \nabla,\Delta_q]\delta d\|_{L^1_t(L^2)}+\|\delta u\cdot \nabla d^1\|_{L^1_t(B_{2,1}^{1/2})}
\nonumber\\
&+\|\nabla d^2\cdot\nabla\delta d\cdot d^2\|_{L^1_t(B_{2,1}^{1/2})}+\|\nabla \delta d\cdot\nabla d^1\cdot d^2\|_{L^1_t(B_{2,1}^{1/2})}+\|\nabla d^1\cdot\nabla d^1\cdot \delta d\|_{L^1_t(B_{2,1}^{1/2})}.
\end{align}
Applying Lemma \ref{daishu1}, \ref{jiaohuanzi} gives that
\begin{align}\label{}
 \sum_{q\in\Z}2^{{q}/2} \|[u^2 \cdot \nabla,\Delta_q]\delta d\|_{L^1_t(L^2)}\le C\int_0^t\|\delta d\|_{B_{2,1}^{1/2}}\|u^2\|_{B_{2,1}^{5/2}}d\tau,
\end{align}
\begin{align}\label{}
\|\delta u\cdot \nabla d^1\|_{L^1_t(B_{2,1}^{1/2})}\le& C\int_0^t\|\delta u\|_{B_{2,1}^{1/2}}\|d^1\|_{B_{2,1}^{5/2}}d\tau\nonumber\\
\le& C\int_0^t\|\delta u\|_{B_{2,1}^{-1/2}}^{1/2}\|\delta u\|_{B_{2,1}^{3/2}}^{1/2}\|d^1\|_{B_{2,1}^{3/2}}^{1/2}\|d^1\|_{B_{2,1}^{7/2}}^{1/2}d\tau\nonumber\\
\le& \varepsilon \|\delta u\|_{L^1_t(B_{2,1}^{3/2})}
+C\int_0^t\|\delta u\|_{B_{2,1}^{-1/2}}\|d^1\|_{B_{2,1}^{3/2}}\|d^1\|_{B_{2,1}^{7/2}}d\tau,
\end{align}
\begin{align}\label{}
\|\nabla d^2\cdot\nabla\delta d\cdot d^2\|_{L^1_t(B_{2,1}^{1/2})}\le& C\int_0^t\|\delta d\|_{B_{2,1}^{3/2}}\|d^2\|_{B_{2,1}^{3/2}}\|d^2\|_{B_{2,1}^{5/2}}d\tau\nonumber\\
\le& C\int_0^t\|\delta d\|_{B_{2,1}^{1/2}}^{1/2}\|\delta d\|_{B_{2,1}^{5/2}}^{1/2}\|d^2\|_{B_{2,1}^{3/2}}^{3/2}\|d^2\|_{B_{2,1}^{7/2}}^{1/2}d\tau\nonumber\\
\le& \varepsilon \|\delta d\|_{L^1_t(B_{2,1}^{5/2})} +C\int_0^t\|\delta d\|_{B_{2,1}^{1/2}}\|d^2\|_{B_{2,1}^{3/2}}^3\|d^2\|_{B_{2,1}^{7/2}}d\tau,
\end{align}
\begin{align}\label{}
\|\nabla \delta d\cdot\nabla d^1\cdot d^2\|_{L^1_t(B_{2,1}^{1/2})}\le& C\int_0^t\|\delta d\|_{B_{2,1}^{3/2}}\|d^2\|_{B_{2,1}^{3/2}}\|d^1\|_{B_{2,1}^{5/2}}d\tau\nonumber\\
\le& C\int_0^t\|\delta d\|_{B_{2,1}^{1/2}}^{1/2}\|\delta d\|_{B_{2,1}^{5/2}}^{1/2}\|d^2\|_{B_{2,1}^{3/2}}\|d^1\|_{B_{2,1}^{3/2}}^{1/2}\|d^1\|_{B_{2,1}^{7/2}}^{1/2}d\tau\nonumber\\
\le& \varepsilon \|\delta d\|_{L^1_t(B_{2,1}^{5/2})} +C\int_0^t\|\delta d\|_{B_{2,1}^{1/2}}\|d^2\|_{B_{2,1}^{3/2}}\|d^1\|_{B_{2,1}^{3/2}}\|d^1\|_{B_{2,1}^{7/2}}d\tau,
\end{align}
\begin{align}\label{}
\|\nabla d^1\cdot\nabla d^1\cdot \delta d\|_{L^1_t(B_{2,1}^{1/2})}\le& C\int_0^t\|\delta d\|_{B_{2,1}^{1/2}}\|d^1\|_{B_{2,1}^{5/2}}^2d\tau\nonumber\\
\le& C\int_0^t\|\delta d\|_{B_{2,1}^{1/2}}\|d^1\|_{B_{2,1}^{3/2}}\|d^1\|_{B_{2,1}^{7/2}}d\tau.
\end{align}
Choosing $\varepsilon$ small enough and inserting the above estimates into (\ref{5.3}),  we have
\begin{align}\label{5.4}
&\|\delta d\|_{\tilde{L}^\infty_t(B_{2,1}^{1/2})}+\|\delta d\|_{L^1_t(B_{2,1}^{5/2})}\nonumber\\
&\lesssim
\int_0^t(\|\delta u\|_{B_{2,1}^{-1/2}}+\|\delta d\|_{B_{2,1}^{1/2}})( \|u^2\|_{B_{2,1}^{5/2}}+\|d^2\|_{B_{2,1}^{3/2}}^3\|d^2\|_{B_{2,1}^{7/2}}+(1+\|d^2\|_{B_{2,1}^{3/2}}) \|d^1\|_{B_{2,1}^{3/2}}\|d^1\|_{B_{2,1}^{7/2}})d\tau.
\end{align}

Combining (\ref{5.12})  with  (\ref{5.4}), applying Gronwall's inequality and using (\ref{5.0}) implies
$\delta a=\delta u = \delta d= 0$ for all $t \in [0, T ].$
\begin{center}
\section{Global well-posedness of Theorem 1.1}
\end{center}
\subsection{The estimate of the transport equation}
In this section, we shall investigate the following transport equation:
\begin{align}\label{zhiliangfangcheng}
\left\{\begin{aligned}
&\partial_ta+u\cdot\nabla a=0,\\
&\mathrm{\mathrm{div}}u=0,\\
&a|_{t=0}=a_0.
\end{aligned}\right.
\end{align}
\begin{proposition}\label{zhiliang}
 Let $u=(u^h,u^3)\in\tilde{L}_T^\infty({\dot{B}_{2,1}^{1/2}}(\R^3))\cap L_T^1({\dot{B}_{2,1}^{5/2}}(\R^3))$ with $\mathrm{div}\,u=0$ and $a_0\in \dot{B}_{2,1}^{3/2}(\R^3)$.
Then $\mathrm{(\ref{zhiliangfangcheng})}$ has a unique solution
$a\in C([0,T ];\dot{B}_{2,1}^{3/2}(\R^3))$ so that
\begin{align}\label{zhiliangguji}
\|a\|_{\tilde{L}_t^{\infty}(\dot{B}_{2,1}^{3/2})}
\le\|a_0\|_{\dot{B}_{2,1}^{3/2}}
+C\|a\|_{\tilde{L}_t^{\infty}(\dot{B}_{2,1}^{3/2})}(\|u^3\|_{{L}_{t}^1({\dot{B}_{2,1}^{5/2}})}^{{1}/{2}}
\|u^h\|_{{L}_{t}^1({\dot{B}_{2,1}^{5/2}})}^{{1}/{2}}+\|u^h\|_{L_t^1({\dot{B}_{2,1}^{5/2}})})
\end{align}
for any $ t\in (0,T ]$.
\end{proposition}
$\mathbf{Proof:}$
The existence and uniqueness of solutions to (\ref{zhiliangfangcheng}) essentially follow from the estimate (\ref{zhiliangguji}) for some appropriate solutions to (\ref{zhiliangfangcheng}). For simplicity, we just present the estimate (\ref{zhiliangguji}) for smooth enough solutions of (\ref{zhiliangfangcheng}) . Thanks to Bony's decomposition (\ref{bony}), we obtain
\begin{align}\label{ZZ1}
u\cdot\nabla a=T_u\cdot \nabla a+\mathcal{R}(u,\nabla a)
\end{align}
Applying $\dot{\Delta}_j$ to (\ref{zhiliangfangcheng}) and taking $L^2$ inner product
with $\dot{\Delta}_ja$, we obtain
\begin{align}\label{ZZ2}
\frac12\frac{d}{dt} \|\Delta_ja\|_{L^2}^2+(\Delta_j(T_u\cdot \nabla a),\dot{\Delta}_ja)+(\Delta_j(\mathcal{R}(u,\nabla a)),\dot{\Delta}_ja).
\end{align}
Using a standard commutator's argument and the $L^2$ energy estimate in  \cite{danchincpde2001}, we have
\begin{align}\label{gujishi}
\|\Delta_ja\|_{L^2}\le& \|\Delta_ja_0\|_{L^2}+C\int_0^t(\sum_{|j-j'|\le5}(\|[\Delta_j,S_{j'-1}]\Delta_{j'}\nabla a\|_{L^2}+\|(S_{j'-1}u-S_{j-1}u)\Delta_j\Delta_{j'}\nabla a\|_{L^2})d\tau\nonumber\\
&+C\int_0^t(\sum_{|j-j'|\le5}(\|\Delta_j\mathcal{R}(u,\nabla a)\|_{L^2}))d\tau.
\end{align}
We first get by applying the classical estimate on commutators and (\ref{guanjian2}) with  $ m=r=\infty$ that
\begin{align*}
&\sum_{|j-j'|\le5}(\|[\Delta_j,S_{j'-1}]\Delta_{j'}\nabla a\|_{L^1_t(L^2)}+\|(S_{j'-1}u-S_{j-1}u)\Delta_j\Delta_{j'}\nabla a\|_{L^1_t(L^2)})\nonumber\\
\le& C\sum_{|j-j'|\le5}(\|S_{j'-1}\nabla u^h\|_{L_t^1(L^\infty)}
+\|S_{j'-1}\nabla u^3\|_{L_t^1(L^\infty)})\|\Delta_{j'}a\|_{L^\infty_t(L^2)}\nonumber\\
&+\sum_{|j-j'|\le5}(\|S_{j'-1}\nabla u^h-S_{j-1}\nabla u^h\|_{L_t^1(L^\infty)}
+\|\|S_{j'-1}\nabla u^3-S_{j-1}\nabla u^3\|_{L_t^1(L^\infty)})\|\Delta_{j'}a\|_{L^\infty_t(L^2)}\nonumber\\
\le& C\sum_{|j-j'|\le5}\sum_{j''\le j'-2}(\|\Delta_{j''}\nabla u^h\|_{L_t^1(L^\infty)}
+2^{j''}\|\Delta_{j''} u^3\|_{L_t^1(L^\infty)})\|\Delta_{j'}a\|_{L^\infty_t(L^2)}\nonumber\\
\le& Cd_j2^{-{3j}/{2}} \|a\|_{\tilde{L}_t^{\infty}(\dot{B}_{2,1}^{3/2})}(\|u^3\|_{{L}_{t}^1({\dot{B}_{2,1}^{5/2}})}^{{1}/{2}}
\|u^h\|_{{L}_{t}^1({\dot{B}_{2,1}^{5/2}})}^{{1}/{2}}+\|u^h\|_{L_t^1({\dot{B}_{2,1}^{5/2}})}).
\end{align*}
Thanks to Lemma \ref{bernstein2} and (\ref{guanjian2}) with  $ m=2, r=\infty,$  we obtain
\begin{align*}
&\|\Delta_j\mathcal{R}(u,\nabla a)\|_{L^1_t(L^2)}\nonumber\\
\le &C\sum_{j'\ge j-N_0}(\|S_{j'+2}\nabla_ha\|_{L_t^\infty(L^\infty)}\|\Delta_{j'}u^h\|_{L^\infty_t(L^2)}
+\|S_{j'+2}\partial_3a\|_{L_t^\infty(L_h^\infty(L^2_v))} \|\Delta_{j'}u^3\|_{L^\infty_t(L^2_h(L_v^\infty))})\nonumber\\
\le &C\sum_{j'\ge j-N_0}\sum_{j''\le j'+1}(2^{(5j''/2)}\|\Delta_{j''}a\|_{L_t^\infty(L^2)}\|\Delta_{j'}u^h\|_{L^\infty_t(L^2)}
+2^{(2j'')}\|\Delta_{j''}a\|_{L_t^\infty(L^2)}\|\Delta_{j'}u^3\|_{L^1_t(L^2_h(L_v^\infty))})\nonumber\\
\le& Cd_j2^{-{3j}/{2}} \|a\|_{\tilde{L}_t^{\infty}(\dot{B}_{2,1}^{3/2})}(\|u^3\|_{{L}_{t}^1({\dot{B}_{2,1}^{5/2}})}^{{1}/{2}}
\|u^h\|_{{L}_{t}^1({\dot{B}_{2,1}^{5/2}})}^{{1}/{2}}+\|u^h\|_{L_t^1({\dot{B}_{2,1}^{5/2}})}).
\end{align*}
Taking above estimates into (\ref{gujishi})  and taking summation for $j\in\Z$, we conclude the proof of (\ref{zhiliangguji}).

\subsection{The estimate of $d$}
Applying $\dot{\Delta}_j$ to (\ref{2}), a standard commutator process gives
\begin{align}\label{Z1}
\partial_t\dot{\Delta}_jd-\Delta\dot{\Delta}_jd=-\dot{\Delta}_j(\dot{T}_u\nabla d)-\dot{\Delta}_j(\mathcal{\dot{R}}(u,\nabla d))+\dot{\Delta}_j(|\nabla d|^2d).
\end{align}
Taking  $L^2$ inner product
with ${\dot{\Delta}}_jd$ to the above equation, we obtain
\begin{align}\label{Z2}
\frac12\frac{d}{dt} \|\dot{\Delta}_jd\|_{L^2}^2+C2^{2j}\|\dot{\Delta}_jd\|_{L^2}^2
\le& C(\dot{\Delta}_j(\dot{T}_u \nabla d),{\dot{\Delta}}_jd)+C(\dot{\Delta}_j(\mathcal{\dot{R}}(u,\nabla d)),{\dot{\Delta}}_jd)\nonumber\\
&+C(\dot{\Delta}_j(|\nabla d|^2d),{\dot{\Delta}}_jd),
\end{align}
where we have used the following fact:

there exists a positive constant $C$ so that
$$
-\int_{{\mathbb R}^3}\Delta\dot{\Delta}_jd{\dot{\Delta}}_jddx\ge C2^{2j}\|\dot{\Delta}_jd\|_{L^2}^2.
$$
Using a standard commutator's argument and the $L^2$ energy estimate in  \cite{danchincpde2001}, we have
\begin{align}\label{Z3}
&\|\dot{\Delta}_jd\|_{L^2}+C2^{2j}\|\dot{\Delta}_jd\|_{L_t^1{(L^2)}} \nonumber\\\le& \|\dot{\Delta}_jd_0\|_{L^2}+C\|\dot{\Delta}_j(|\nabla d|^2d)\|_{L_t^1{(L^2)}}+C\|\dot{\Delta}_j(d\times \Delta d)\|_{L_t^1{(L^2)}}\nonumber\\
&+C\int_0^t(\sum_{|j-j'|\le5}(\|[\dot{\Delta}_j,\dot{S}_{j'-1}u]\dot{\Delta}_{j'}\nabla d\|_{L^2}+\|(\dot{S}_{j'-1}u-\dot{S}_{j-1}u)\dot{\Delta}_j\dot{\Delta}_{j'}\nabla d\|_{L^2})d\tau\nonumber\\
&+C\int_0^t(\sum_{|j-j'|\le5}(\|\dot{\Delta}_j\mathcal{\dot{R}}(u,\nabla d)\|_{L^2}))d\tau.
\end{align}
We first get by applying the classical estimate on commutators and (\ref{guanjian2}) with  $ q=r=\infty$ that
\begin{align*}
&\sum_{|j-j'|\le5}(\|[\dot{\Delta}_j,\dot{S}_{j'-1}u]\dot{\Delta}_{j'}\nabla d\|_{L^1_t(L^2)}+\|(\dot{S}_{j'-1}u-\dot{S}_{j-1}u)\dot{\Delta}_j\dot{\Delta}_{j'}\nabla d\|_{L^1_t(L^2)})\nonumber\\
\le& C\sum_{|j-j'|\le5}(\|\dot{S}_{j'-1}\nabla u^h\|_{L_t^1(L^\infty)}
\nonumber\\
&+C\sum_{|j-j'|\le5}(\|\dot{S}_{j'-1}\nabla u^h-\dot{S}_{j-1}\nabla u^h\|_{L_t^1(L^\infty)}
+\|\dot{S}_{j'-1}\nabla u^3-\dot{S}_{j-1}\nabla u^3\|_{L_t^1(L^\infty)})\|\dot{\Delta}_{j'}d\|_{L^\infty_t(L^2)}\nonumber\\
\le& C\sum_{|j-j'|\le5}\sum_{j''\le j'-2}(2^{j''}\|\dot{\Delta}_{j''} u^h\|_{L_t^1(L^\infty)}
+2^{j''}\|\dot{\Delta}_{j''} u^3\|_{L_t^1(L^\infty)})\|\dot{\Delta}_{j'}d\|_{L^\infty_t(L^2)}\nonumber\\
\le& Cd_j2^{-3j/2} \|d\|_{\tilde{L}_t^{\infty}(\dot{B}_{2,1}^{3/2})}(\|u^3\|_{{L}_{t}^1({\dot{B}_{2,1}^{5/2}})}^{1/2}
\|u^h\|_{{L}_{t}^1({\dot{B}_{2,1}^{5/2}})}^{1/2}+\|u^h\|_{L_t^1({\dot{B}_{2,1}^{5/2}})}).
\end{align*}
Thanks to Lemma \ref{bernstein2}, \ref{daishu1} and (\ref{guanjian2}) with  $ q=2, r=\infty,$  we obtain
\begin{align*}
&\|\dot{\Delta}_j\mathcal{\dot{R}}(u,\nabla d)\|_{L^1_t(L^2)}\nonumber\\
\le &C\sum_{j'\ge j-N_0}\Big(\|\dot{S}_{j'+2}\nabla_hd\|_{L_t^\infty(L^\infty)}\|\dot{\Delta}_{j'}u^h\|_{L^\infty_t(L^2)}
\nonumber\\
&+\|\dot{S}_{j'+2}\partial_3d\|_{L_t^\infty(L_h^\infty(L^2_v))} \|\dot{\Delta}_{j'}u^3\|_{L^\infty_t(L^2_h(L_v^\infty))}\Big)\nonumber\\
\le &C\sum_{j'\ge j-N_0}\sum_{j''\le j'+1}\Big(2^{5j''/2}\|\dot{\Delta}_{j''}d\|_{L_t^\infty(L^2)}\|\dot{\Delta}_{j'}u^h\|_{L^\infty_t(L^2)}
\nonumber\\
&+2^{2j''}\|\dot{\Delta}_{j''}d\|_{L_t^\infty(L^2)}\|\dot{\Delta}_{j'}u^3\|_{L^1_t(L^2_h(L_v^\infty))}\Big)\nonumber\\
\le& Cd_j2^{-3j/2} \|d\|_{\tilde{L}_t^{\infty}(\dot{B}_{2,1}^{3/2})}(\|u^3\|_{{L}_{t}^1({\dot{B}_{2,1}^{5/2}})}^{1/2}
\|u^h\|_{{L}_{t}^1({\dot{B}_{2,1}^{5/2}})}^{1/2}+\|u^h\|_{L_t^1({\dot{B}_{2,1}^{5/2}})})
\end{align*}
and
\begin{align}\label{Z5}
\|\dot{\Delta}_j(|\nabla d|^2d)\|_{L_t^1{(L^2)}}\le Cd_j2^{-3j/2}\|d\|_{\tilde{L}^\infty_t(\dot{B}_{2,1}^{3/2})}^2\|d\|_{L^1_t(\dot{B}_{2,1}^{7/2})},
\end{align}

Substituting above estimates into (\ref{gujishi})  and taking summation for $j\in{\mathbb Z}$, we can get
\begin{align}\label{Z6}
\|d\|_{\tilde{L}^\infty_t(\dot{B}_{2,1}^{3/2})}+C\|d\|_{L^1_t(\dot{B}_{2,1}^{7/2})}
\le&\|d_0\|_{\dot{B}_{2,1}^{3/2}}+C\|d\|_{\tilde{L}^\infty_t(\dot{B}_{2,1}^{3/2})}^2\|d\|_{L^1_t(\dot{B}_{2,1}^{7/2})}
\nonumber\\
&+C\|d\|_{\tilde{L}^\infty_t(\dot{B}_{2,1}^{3/2})}(\|u^3\|_{{L}_{t}^1({\dot{B}_{2,1}^{5/2}})}^{1/2}
\|u^h\|_{{L}_{t}^1({\dot{B}_{2,1}^{5/2}})}^{1/2}+\|u^h\|_{L_t^1({\dot{B}_{2,1}^{5/2}})}).\quad\square
\end{align}

\subsection{The estimate of the pressure}
As is well known, the main difficulty in the study of the
well-posedness of incompressible MHD system is to the derive
the estimate for the pressure term.
We first get by taking div
to the second equation of (\ref{1.2}) that
\begin{align}\label{yalifangcheng}
-\Delta P=&\mathrm{div}(a\nabla P)+\mathrm{div}(u\cdot\nabla u)
-\mathrm{div}(a\Delta u)-\mathrm{div}[(1+a)(\nabla d \cdot \Delta d)]\nonumber\\
=& \mathrm{div}(a\nabla P)+ \mathrm{div}_h\mathrm{div}_h(u^h\otimes u^h)+2\partial_3\mathrm{div}_h(u^3u^h)-2\partial_3(u^3\mathrm{div}_hu^h)-\mathrm{div}_h(a\Delta u^h)
\nonumber\\
&-\partial_3(a\Delta u^3)-\mathrm{div}[(1+a)\nabla\cdot (\nabla d \odot \nabla d)].
\end{align}
In what follows, we will give the estimate of the pressure which will be used in the estimates of $(u^h,u^3)$ .

\begin{proposition}\label{yaliguji}
 Let  $a\in\tilde{L}_T^\infty({\dot{B}_{2,1}^{3/2}}(\R^3))$ $u=(u^h,u^3)\in\tilde{L}_T^\infty({\dot{B}_{2,1}^{1/2}}(\R^3))\cap L_T^1({\dot{B}_{2,1}^{5/2}}(\R^3))$, $d\in\tilde{L}_T^\infty({\dot{B}_{2,1}^{3/2}}(\R^3))\cap L_T^1({\dot{B}_{2,1}^{7/2}}(\R^3))$.
Then $\mathrm{(\ref{yalifangcheng})}$ has a unique solution
$\nabla P\in L_T^1({\dot{B}_{2,1}^{1/2}}(\R^3))$ which decays to zero when $|x|\rightarrow\infty$ so that for all $t\in[0,T]$,  there holds
\begin{align}\label{yali}
\|\nabla P\|_{L^1_t(\dot{B}_{2,1}^{1/2})}
\le& C Y(a,u,d),
\end{align}
provided that  $C\|a\|_{\tilde{L}_T^{\infty}(\dot{B}_{2,1}^{3/2})}\le1/2,$
where $Y(a,u,d)$ is defined by
\begin{align*}
Y(a,u,d)
\triangleq& C\Bigg(
\|u^h\|_{\tilde{L}^\infty_t(\dot{B}_{2,1}^{1/2})}
\|u^h\|_{L^1_t(\dot{B}_{2,1}^{5/2})}+\|a\|_{\tilde{L}_t^{\infty}(\dot{B}_{2,1}^{3/2})}(\|u^3\|_{L^1_t(\dot{B}_{2,1}^{5/2})}+
\|u^h\|_{L^1_t(\dot{B}_{2,1}^{5/2})})\nonumber\\
&+F(u^3,u^h)+(1+\|a\|_{\tilde{L}_t^{\infty}(\dot{B}_{2,1}^{3/2})})
\|d\|_{\tilde{L}^\infty_t(\dot{B}_{2,1}^{3/2})}
\|d\|_{L^1_t(\dot{B}_{2,1}^{7/2})}\Bigg).
\end{align*}
\end{proposition}
$\bf{Proof\,\, of\,\, Proposition \,\ref{yaliguji}.}$  As both the existence and uniqueness parts of Proposition \ref{yaliguji} basically follows from the uniform estimate (\ref{yali}) for appropriate approximate solutions of (\ref{yalifangcheng}) . For simplicity, we just prove (\ref{yali}) for smooth enough solutions of (\ref{yalifangcheng}).
Firstly,
applying $\dot{\Delta}_j$ to (\ref{yalifangcheng}) and using Lemma \ref{bernstein2} we have
\begin{align}\label{yaligujishi}
\|\dot{\Delta}_j(\nabla P)\|_{L^1_t(L^2)}
\lesssim& \|\dot{\Delta}_j(a\nabla P)\|_{L^1_t(L^2)}+
 2^j\|\dot{\Delta}_j(u^h\otimes u^h)\|_{L^1_t(L^2)}
 + 2^j\|\dot{\Delta}_j(u^3u^h)\|_{L^1_t(L^2)}
 \nonumber\\
&+\|\dot{\Delta}_j(u^3\mathrm{div}_hu^h)\|_{L^1_t(L^2)}+\|\dot{\Delta}_j(a\Delta u^h)\|_{L^1_t(L^2)} \nonumber\\
&+\|\dot{\Delta}_j(a\Delta u^3)\|_{L^1_t(L^2)}
+\|\dot{\Delta}_j[(1+a)\nabla\cdot (\nabla d \odot \nabla d)]\|_{L^1_t(L^2)}.
\end{align}
Applying  Lemma \ref{daishu1} and Lemma \ref{zhongyaoguji} gives
\begin{align}\label{P2}
\|\dot{\Delta}_j(a\nabla P)\|_{L^1_t(L^2)}\lesssim d_j2^{(-j/2)} \|a\|_{\tilde{L}_t^{\infty}(\dot{B}_{2,1}^{3/2})}\|\nabla P\|_{L^1_t(\dot{B}_{2,1}^{1/2})},
\end{align}
\begin{align}\label{P3}
2^j\|\dot{\Delta}_j(u^h\otimes u^h)\|_{L^1_t(L^2)}
\lesssim d_j2^{(-j/2)} (\|u^h\|_{\tilde{L}^\infty_t(\dot{B}_{2,1}^{1/2})}
\|u^h\|_{L^1_t(\dot{B}_{2,1}^{5/2})}),
\end{align}
\begin{align}\label{P4}
&2^j\|\dot{\Delta}_j(u^3u^h)\|_{L^1_t(L^2)}
 +\|\dot{\Delta}_j(u^3\mathrm{div}_hu^h)\|_{L^1_t(L^2)}
\lesssim d_j2^{(-j/2)}F(u^3,u^h),
\end{align}
\begin{align}\label{P5}
\|\dot{\Delta}_j(a\Delta u^h)\|_{L^1_t(L^2)}+\|\dot{\Delta}_j(a\Delta u^3)\|_{L^1_t(L^2)}
\lesssim d_j2^{(-j/2)}\|a\|_{\tilde{L}_t^{\infty}(\dot{B}_{2,1}^{3/2})}(\|u^h\|_{L^1_t(\dot{B}_{2,1}^{5/2})}
+\|u^3\|_{L^1_t(\dot{B}_{2,1}^{5/2})}),
\end{align}
\begin{align}\label{P6}
\|\dot{\Delta}_j[(1+a)\nabla\cdot (\nabla d \odot \nabla d)]\|_{L^1_t(L^2)} \lesssim&d_j2^{(-j/2)}(1+\|a\|_{\tilde{L}_t^{\infty}(\dot{B}_{2,1}^{3/2})})\|d\|_{\tilde{L}^\infty_t(\dot{B}_{2,1}^{3/2})}
\|d\|_{L^1_t(\dot{B}_{2,1}^{7/2})}.
\end{align}

Taking above estimates (\ref{P2}) - (\ref{P6}) into (\ref{yaligujishi}), we have
\begin{align}\label{P10}
\|\nabla P\|_{L^1_t(\dot{B}_{2,1}^{1/2})}
\le& C\Bigg(\|a\|_{\tilde{L}_t^{\infty}(\dot{B}_{2,1}^{3/2})}\|\nabla P\|_{L^1_t(\dot{B}_{2,1}^{1/2})}
+\|u^h\|_{\tilde{L}^\infty_t(\dot{B}_{2,1}^{1/2})}
\|u^h\|_{L^1_t(\dot{B}_{2,1}^{5/2})}+F(u^3,u^h)\nonumber\\
&
+\|a\|_{\tilde{L}_t^{\infty}(\dot{B}_{2,1}^{3/2})}(\|u^3\|_{L^1_t(\dot{B}_{2,1}^{5/2})}+
\|u^h\|_{L^1_t(\dot{B}_{2,1}^{5/2})})
+(1+\|a\|_{\tilde{L}_t^{\infty}(\dot{B}_{2,1}^{3/2})})
\|d\|_{\tilde{L}^\infty_t(\dot{B}_{2,1}^{3/2})}
\|d\|_{L^1_t(\dot{B}_{2,1}^{7/2})}
\Bigg).
\end{align}
So provided that $C\|a\|_{\tilde{L}_t^{\infty}(\dot{B}_{2,1}^{3/2})}\le1/2,$
we can conclude the proof of (\ref{yali}).
\subsection{The estimates of the horizontal component.}
Thanks to the second equation of (\ref{1.2}) again, we have
\begin{align}\label{shuipingfangxiang}
\partial_{t} u^h
+u \cdot \nabla u^h-(1+a)(\Delta u^h-\nabla_hP)=&-(1+a)\nabla\cdot (\nabla d \odot \nabla_h d).
\end{align}
Applying $\dot{\Delta}_j$ to (\ref{shuipingfangxiang}) and taking $L^2$ inner product
with $\dot{\Delta}_ju^h$, we obtain
\begin{align}\label{a11}
&\frac{d}{dt}\|\dot{\Delta}_ju^h\|_{L^2}
+C_12^{2j}\|\dot{\Delta}_ju^h\|_{L^2}\nonumber\\
\lesssim&\|\dot{\Delta}_ju_0^h\|_{L^2}+\|\dot{\Delta}_j(u \cdot \nabla u^h)\|_{L^2}+\|\dot{\Delta}_j((1+a)\nabla\cdot (\nabla d \odot \nabla_h d))\|_{L^2}\nonumber\\
&+\|\dot{\Delta}_j(a\Delta u^h)\|_{L^2}+\|\dot{\Delta}_j((1+a)\nabla_h
P)\|_{L^2}.
\end{align}
Integrating the resulting inequality over [0,t]and using the fact $\mathrm{\mathrm{div}}u=0,$  one shows that
\begin{align}\label{a13}
&\|\dot{\Delta}_ju^h\|_{L_t^\infty(L^2)}
+ 2^{2j}(\|\dot{\Delta}_ju^h\|_{L^1_t(L^2)}
\nonumber\\
\lesssim&\|\dot{\Delta}_ju_0^h\|_{L^2}+\|\dot{\Delta}_j(u \cdot \nabla u^h)\|_{L^1_t(L^2)}
+\|\dot{\Delta}_j((1+a)\nabla\cdot (\nabla d \odot \nabla_h d))\|_{L^1_t(L^2)}
\nonumber\\
&
+\|\dot{\Delta}_j((1+a)\nabla_hP)\|_{L^1_t(L^2)}
+\|\dot{\Delta}_j(a\Delta u^h)\|_{L^1_t(L^2)}\nonumber\\
\lesssim&\|\dot{\Delta}_ju_0^h\|_{L^2}+2^j\|u^h\otimes u^h\|_{L^1_t(L^2)}+2^j\|u^3 u^h\|_{L^1_t(L^2)}
+\|\dot{\Delta}_j((1+a)\nabla\cdot (\nabla d \odot \nabla_h d))\|_{L^1_t(L^2)}
\nonumber\\
&
+\|\dot{\Delta}_j((1+a)\nabla_hP)\|_{L^1_t(L^2)}
+\|\dot{\Delta}_j(a\Delta u^h)\|_{L^1_t(L^2)}.
\end{align}
Therefore, we can get from Proposition \ref{yaliguji} and the estimates (\ref{P3})-(\ref{P6})  that
 \begin{align}\label{a18}
\|u^h\|_{\tilde{L}^\infty_t(\dot{B}_{2,1}^{1/2})}
+\|u^h\|_{L^1_t(\dot{B}_{2,1}^{5/2})}
\le \|u_0^h\|_{\dot{B}_{2,1}^{1/2}}+Y(a,u,d).
\end{align}
provided that $C\|a\|_{\tilde{L}_t^{\infty}(\dot{B}_{2,1}^{3/2})}\le 1/2$.
\subsection{The estimates of the vertical component}
In this section, we will use the fact that the velocity field
equation on the vertical component is a linear equation with
 coefficients depending  only on the horizontal components and $a$.
 Thanks to $(\ref{1.1})_2$, and $\mathrm{\mathrm{div}}u=0$,
 we only consider the vertical component of the system to get
\begin{align}\label{V1}
\partial_{t} u^3
+u \cdot \nabla u^3-(1+a)(\Delta u^3-\partial_3P)=&-(1+a)\nabla\cdot (\nabla d \odot \partial_3 d).
\end{align}
Applying the operator $\dot{\Delta}_j$ to above equations and taking $L^2$ inner product of the resulting
equation with $\dot{\Delta}_ju^3$ ,
we get by a similar derivation of
(\ref{a13}) that
\begin{align}\label{V2}
&\|\dot{\Delta}_ju^3\|_{L_t^\infty(L^2)}
+2^{2j}\|\dot{\Delta}_ju^3\|_{L_t^1(L^2)}
\nonumber\\
\lesssim&\|\dot{\Delta}_ju_0^3\|_{L^2}
+\|\dot{\Delta}_j(u \cdot \nabla u^3)\|_{L_t^1(L^2)}+\|\dot{\Delta}_j((1+a)\partial_3P)\|_{L_t^1(L^2)}\nonumber\\
&+\|\dot{\Delta}_j((1+a)\nabla\cdot (\nabla d \odot \partial_3 d))\|_{L_t^1(L^2)}
+\|\dot{\Delta}_j(a\Delta u^3)\|_{L_t^1(L^2)}.
\end{align}
Applying Lemma \ref{daishu1} and Young's inequality gives
\begin{align}\label{V3}
\|\dot{\Delta}_j(u \cdot \nabla u^3)\|_{L_t^1(L^2)}
\le& C 2^j\|\dot{\Delta}_j(u^hu^3)\|_{L_t^1(L^2)}\nonumber\\
\le& Cd_j2^{(-j/2)}(\|u^3\|_{\tilde{L}_t^2({\dot{B}_{2,1}^{3/2}})}\|u^h\|_{\tilde{L}_t^2({\dot{B}_{2,1}^{3/2}})}
)\nonumber\\
\le
&Cd_j2^{(-j/2)}(\|u^3\|_{\tilde{L}^\infty_t(\dot{B}_{2,1}^{1/2})}
+
\|u^3\|_{L^1_t(\dot{B}_{2,1}^{5/2})})(\|u^h\|_{\tilde{L}^\infty_t(\dot{B}_{2,1}^{1/2})}
+
\|u^h\|_{L^1_t(\dot{B}_{2,1}^{5/2})}),
\end{align}
\begin{align}\label{V4}
\|\dot{\Delta}_j(a\Delta u^3)\|_{L_t^1(L^2)}\le& Cd_j2^{(-j/2)}\|a\|_{\tilde{L}_t^{\infty}(\dot{B}_{2,1}^{3/2})}\|u^3\|_{L^1_t(\dot{B}_{2,1}^{5/2})},\nonumber\\
\|\dot{\Delta}_j((1+a)\partial_3P)\|_{L_t^1(L^2)}\le& Cd_j2^{(-j/2)} \|a\|_{\tilde{L}_t^{\infty}(\dot{B}_{2,1}^{3/2})}Y(a,u,b)
\end{align}
and
\begin{align}\label{V5}
&\|\dot{\Delta}_j((1+a)\nabla\cdot (\nabla d \odot \partial_3 d))\|_{L_t^1(L^2)}\nonumber\\
\le&Cd_j2^{(-j/2)}
(1+\|a\|_{\tilde{L}_t^{\infty}(\dot{B}_{2,1}^{3/2})})
(\|d\|_{\tilde{L}^\infty_t\dot{B}_{2,1}^{3/2})}
\|d\|_{L^1_t(\dot{B}_{2,1}^{7/2})}.
\end{align}

Inserting estimates (\ref{V3}) - (\ref{V5}) into (\ref{V2}), we have
\begin{align}\label{chuizhifangxiang}
\|u^3\|_{\tilde{L}^\infty_t(\dot{B}_{2,1}^{1/2})}
+\|u^3\|_{L^1_t(\dot{B}_{2,1}^{5/2})}
\le& \|u_0^3\|_{\dot{B}_{2,1}^{1/2}}
+Y(a,u,b)\nonumber\\
&+C(\|u^3\|_{\tilde{L}^\infty_t(\dot{B}_{2,1}^{1/2})}+
\|u^3\|_{L^1_t(\dot{B}_{2,1}^{5/2})})(\|u^h\|_{\tilde{L}^\infty_t(\dot{B}_{2,1}^{1/2})}+
\|u^h\|_{L^1_t(\dot{B}_{2,1}^{5/2})}).
\end{align}
provided that $C\|a\|_{\tilde{L}_t^{\infty}(\dot{B}_{2,1}^{3/2})}\le 1/2$.

\subsection{The proof of Theorem 1.1}
The goal of this section is to present the proof of Theorem 1.1. In fact, given $a_0\in \dot{B}_{2,1}^{3/2}(\R^3)$, $u_0\in \dot{B}_{2,1}^{1/2}(\R^3)$, $d_0\in \dot{B}_{2,1}^{3/2}(\R^3)$ with $\|a_0\|_{\dot{B}_{2,1}^{3/2}}$
being sufficiently small, it follows by a similar
argument in section 4  that there exists a positive time $T$ so that (1.2) has a local solution $(a,u,d)$ with
\begin{align}\label{QQ1}
&a\in C([0,T];\dot{B}_{2,1}^{3/2} )\cap
\tilde{L}^\infty([0,T];\dot{B}_{2,1}^{3/2}), \nonumber\\
&u\in C([0,T];\dot{B}_{2,1}^{1/2})\cap
\tilde{L}^\infty([0,T];\dot{B}_{2,1}^{1/2})\cap
L^1([0,T];\dot{B}_{2,1}^{5/2}),\nonumber\\
&d\in C([0,T];\dot{B}_{2,1}^{3/2})\cap
\tilde{L}^\infty([0,T];\dot{B}_{2,1}^{3/2})\cap
L^1([0,T];\dot{B}_{2,1}^{7/2}).
\end{align}
We denote $T^\ast$ to be the largest possible time such that there holds (4.1).
Then, the proof of Theorem 1.1 is reduced to show that $T^\ast=\infty$ under the assumptions of (\ref{1.2}), (\ref{1.3}).
Let $\epsilon,\eta$ be  small enough positive constants, we define $T^{\ast\ast}$ by
\begin{align}\label{QQ2}
T^{\ast\ast} =\sup\Bigg\{t\in[0,T^\ast):
\|a\|_{\tilde{L}_t^{\infty}(\dot{B}_{2,1}^{3/2})}&\le 4\|a_0\|_{\dot{B}_{2,1}^{3/2}},\nonumber\\
\|u^3\|_{\tilde{L}^\infty_t(\dot{B}_{2,1}^{{1/2}})}
+\|u^3\|_{L^1_t(\dot{B}_{2,1}^{{5/2}})}
&\le 4 (\|u_0^3\|_{\dot{B}_{2,1}^{1/2}}+\epsilon)\nonumber\\
\|u^h\|_{\tilde{L}^\infty_t(\dot{B}_{2,1}^{{1/2}})}+\|d\|_{\tilde{L}^\infty_t(\dot{B}_{2,1}^{{3/2}})}
+\|u^h\|_{L^1_t(\dot{B}_{2,1}^{{5/2}})}+\|d\|_{L^1_t(\dot{B}_{2,1}^{{5/2}})}
&
\le 4 (\|u_0^h\|_{\dot{B}_{2,1}^{1/2}}+\|d\|_{\dot{B}_{2,1}^{3/2}})+\eta\triangleq H_0
\Bigg\}.
\end{align}
In what follows, we will prove that $T^{\ast\ast}=T^\ast$ under the assumptions of (\ref{1.3}), (\ref{1.4}).

If not, we assume that $T^{\ast\ast}<T^\ast,$ and for $\forall t\le T^{\ast\ast}$,
we get from (\ref{zhiliangguji}) that
\begin{align}\label{QQ3}
\|a\|_{\tilde{L}_t^{\infty}(\dot{B}_{2,1}^{3/2})}
\le\|a_0\|_{\dot{B}_{2,1}^{3/2}}
+C\|a\|_{\tilde{L}_t^{\infty}(\dot{B}_{2,1}^{3/2})}(H_0+(\|u_0^3\|_{\dot{B}_{2,1}^{1/2}}+\epsilon)^{1/2}H_0^{1/2}).
\end{align}
By taking
\begin{align}\label{QQ4}
C(H_0+(\|u_0^3\|_{\dot{B}_{2,1}^{1/2}}+\epsilon)^{1/2}H_0^{1/2})<1/2,
\end{align}
we have
\begin{align}\label{QQ5}
\|a\|_{\tilde{L}_t^{\infty}(\dot{B}_{2,1}^{3/2})}
\le2\|a_0\|_{\dot{B}_{2,1}^{3/2}}\hspace{0.5cm} \forall t\le T^{\ast\ast}.
\end{align}
Thanks to (\ref{chuizhifangxiang}) and (\ref{QQ2}), we have
\begin{align}\label{QQ6}
\|u^3\|_{\tilde{L}^\infty_t(\dot{B}_{2,1}^{1/2})}
+\|u^3\|_{L^1_t(\dot{B}_{2,1}^{5/2})}
\le& \|u_0^3\|_{\dot{B}_{2,1}^{1/2}}
+C\|a_0\|_{\dot{B}_{2,1}^{3/2}}H_0+C\|a_0\|_{\dot{B}_{2,1}^{3/2}}\|u^3\|_{L^1_t(\dot{B}_{2,1}^{5/2})}
\nonumber\\
&+C(1+\|a_0\|_{\dot{B}_{2,1}^{3/2}})H_0^2
+C(\|u_0^3\|_{\dot{B}_{2,1}^{1/2}}+\epsilon)^{1/2}H_0^{3/2},
\end{align}
where the following estimate has been used:
\begin{align}\label{QQ7}
F(u^3,u^h)
=&\|u^3\|_{{L}_{t}^1({\dot{B}_{2,1}^{5/2}})}^{1/2}\|u^h\|_{{L}_{t}^1({\dot{B}_{2,1}^{5/2}})}^{3/2}\nonumber\\
&+(\|u^3\|_{{L}_{t}^1({\dot{B}_{2,1}^{5/2}})}+\|u^3\|_{\tilde{L}_t^\infty({\dot{B}_{2,1}^{1/2}})})^{1/2}
(\|u^h\|_{{L}_{t}^1({\dot{B}_{2,1}^{5/2}})}+\|u^h\|_{\tilde{L}_t^\infty({\dot{B}_{2,1}^{1/2}})})^{3/2}\nonumber\\
\le&C(\|u_0^3\|_{\dot{B}_{2,1}^{1/2}}+\epsilon)^{1/2}H_0^{3/2}.
\end{align}
One hand,
for $\forall t\le T^{\ast\ast}$, while taking
\begin{align}\label{QQ8}
&C\|a_0\|_{\dot{B}_{2,1}^{3/2}}H_0\le1/8,\quad\quad
C\|a_0\|_{\dot{B}_{2,1}^{3/2}}\le1/8,\quad\quad C(1+\|a_0\|_{\dot{B}_{2,1}^{3/2}})H_0^2\le1/8,\nonumber\\ &C(\|u_0^3\|_{\dot{B}_{2,1}^{1/2}}+\epsilon)^{1/2}H_0^{3/2}\le1/8.
\end{align}
We can deduce from (\ref{QQ6}) that
\begin{align}\label{QQ12}
\|u^3\|_{\tilde{L}^\infty_t(\dot{B}_{2,1}^{1/2})}
+\|u^3\|_{L^1_t(\dot{B}_{2,1}^{5/2})}
\le8/7(\|u_0^3\|_{\dot{B}_{2,1}^{1/2}}+\epsilon).
\end{align}
On the other hand, for $t\le T^{\ast\ast}$, we get from (\ref{Z6}), (\ref{a18}) that
\begin{align}\label{QQ13}
&\|u^h\|_{\tilde{L}^\infty_t(\dot{B}_{2,1}^{1/2})}
+\|u^h\|_{L^1_t(\dot{B}_{2,1}^{5/2})}+\|d\|_{\tilde{L}^\infty_t(\dot{B}_{2,1}^{3/2})}
+\|d\|_{L^1_t(\dot{B}_{2,1}^{7/2})}
\nonumber\\
\le& \|u_0^h\|_{\dot{B}_{2,1}^{1/2}}
+\|d_0\|_{\dot{B}_{2,1}^{3/2}}
+CH_0\|u^h\|_{L^1_t(\dot{B}_{2,1}^{5/2})}+C\|a_0\|_{\dot{B}_{2,1}^{3/2}}\|u^h\|_{L^1_t(\dot{B}_{2,1}^{5/2})}\nonumber\\
&+C\|a_0\|_{\dot{B}_{2,1}^{3/2}}(\|u_0^3\|_{\dot{B}_{2,1}^{1/2}}+\epsilon)+CH_0^3+C(1+\|a_0\|_{\dot{B}_{2,1}^{3/2}})H_0^2
+C(\|u_0^3\|_{\dot{B}_{2,1}^{1/2}}+\epsilon)^{1/2}H_0^{3/2}.
\end{align}

Choosing
\begin{align}\label{QQ15}
&CH_0(1+H_0^2)\le1/8,\quad\quad C\|a_0\|_{\dot{B}_{2,1}^{3/2}}\le1/8,\quad\quad
C\|a_0\|_{\dot{B}_{2,1}^{3/2}}(\|u_0^3\|_{\dot{B}_{2,1}^{1/2}}+\epsilon)\le1/8,\nonumber\\
&C(1+\|a_0\|_{\dot{B}_{2,1}^{3/2}})H_0^2\le1/8,\quad\quad C(\|u_0^3\|_{\dot{B}_{2,1}^{1/2}}+\epsilon)^{1/2}H_0^{3/2}\le1/8.
\end{align}
we get from (\ref{QQ13}) that
\begin{align}\label{QQ16}
\|u^h\|_{\tilde{L}^\infty_t(\dot{B}_{2,1}^{1/2})}
+\|u^h\|_{L^1_t(\dot{B}_{2,1}^{5/2})}+\|d\|_{\tilde{L}^\infty_t(\dot{B}_{2,1}^{3/2})}
+\|d\|_{L^1_t(\dot{B}_{2,1}^{7/2})}
\le4H_0/3.
\end{align}
Combining with (\ref{QQ4}), (\ref{QQ8}) and (\ref{QQ15}), we can get (\ref{QQ5}), (\ref{QQ12}) and (\ref{QQ16}) hold if we take C large enough in (\ref{1.2}), (\ref{1.3}),  this contradicts with the definition  (\ref{QQ2}), thus we conclude that $T^{\ast\ast}=T^{\ast}$.
Consequently, we complete the proof of Theorem 1.1 by standard continuation argument.

\end{document}